\documentclass[11pt,a4paper]{article}
\usepackage[margin=27mm]{geometry}
\usepackage[T1]{fontenc}
\usepackage[utf8]{inputenc}
\usepackage{lmodern,amsmath,amssymb,amsthm,mathtools}
\usepackage[expansion=false]{microtype}
\usepackage{xcolor}
\definecolor{linkblue}{RGB}{25,61,127}
\usepackage[unicode,bookmarksopen=true,bookmarksnumbered=true,
 colorlinks=true,linkcolor=linkblue,citecolor=linkblue,urlcolor=linkblue,
 linktoc=all,pdfborder={0 0 0}]{hyperref}
\usepackage{bookmark}
\hypersetup{pdftitle={Bernstein--Sato polynomials of generic Pfaffian ideals},
 pdfauthor={Yifan Chen and Huaiqing Zuo},
 pdfsubject={Roots, multiplicities, and strong motivic monodromy},
 pdfkeywords={Pfaffian ideal, Bernstein-Sato polynomial, motivic zeta function}}
\newtheorem{maintheorem}{Theorem}

\newtheorem{theorem}{Theorem}[section]
\newtheorem{proposition}[theorem]{Proposition}
\newtheorem{lemma}[theorem]{Lemma}
\newtheorem{corollary}[theorem]{Corollary}
\theoremstyle{definition}
\newtheorem{definition}[theorem]{Definition}
\newtheorem{remark}[theorem]{Remark}
\newcommand{\C}{\mathbb C}\newcommand{\Z}{\mathbb Z}\newcommand{\N}{\mathbb N}
\newcommand{\D}{\mathcal D}\newcommand{\Pf}{\operatorname{Pf}}\newcommand{\GL}{\operatorname{GL}}

\newcommand{\pr}{\mathfrak p}\newcommand{\mult}{\operatorname{mult}}
\DeclareMathOperator{\tr}{tr}\DeclareMathOperator{\sgn}{sgn}
\numberwithin{equation}{section}
\title{Bernstein--Sato polynomials of generic Pfaffian ideals}
\author{Yifan Chen and Huaiqing Zuo}
\date{}
\begin{document}
\maketitle
\begin{abstract}
We determine all roots and their multiplicities in the  Bernstein--Sato polynomial of every generic Pfaffian ideal. In the nonmaximal case, the roots lie in explicit rational
windows indexed by the rank strata. Integral roots retain their label
multiplicities, whereas every nonintegral root is simple, including
values with several labels. The proof combines graph-module relations,
skew Capelli operators, central characters, and a weighted Mellin
integral. We deduce the strong motivic monodromy conjecture for
Pfaffian ideals, with pole orders bounded by Bernstein--Sato root
multiplicities. For Pfaffians of size $2r$, the largest possible root
multiplicity as the matrix size varies is
$\lfloor(\sqrt{8r+1}-1)/2\rfloor$. In the nonmaximal case we also prove
that the nonunipotent part of the Verdier specialization has zero
nilpotent monodromy.

Keywords. Bernstein-Sato polynomial, Pfaffian ideal, root multiplicities, strong monodromy conjecture.

MSC(2020). 14F10, 14M12, 14E18, 32C38.
\end{abstract}

\begingroup\small
\tableofcontents
\endgroup

\section{Introduction}\label{sec:introduction}
The Bernstein--Sato polynomial records the parameters at which a
functional equation for complex powers of a polynomial degenerates.
Its extension to arbitrary ideals by Budur, Musta\c t\u a, and Saito
\cite{BMS} connects differential operators with multiplier ideals and
Verdier specialization. For ideals defined by rank conditions, the
symmetry and the transverse rank stratification make several geometric
invariants explicit. Determining the complete Bernstein--Sato
polynomial nevertheless requires control of repeated roots, which is
not supplied by the set of monodromy eigenvalues alone.

We study the ordinary ideal $\pr_r=\pr_r^{(N)}$ generated by the
$2r\times2r$ principal Pfaffians of a generic complex alternating
$N\times N$ matrix $X$. Set
\[
 m=\lfloor N/2\rfloor,\qquad \epsilon=N-2m\in\{0,1\},\qquad
 D_0=\binom N2,
 \qquad 1\le r\le m.
\]
Thus $V(\pr_r)$ is the locus of matrices of rank less than $2r$.
Our polynomial $b_{\pr_r}(s)$ is monic and uses the unshifted ideal
convention of \cite{BMS}, recalled in Section~\ref{sec:graph}.

For $r=m$ and $N=2m$, the ideal is principal and the classical
Pfaffian formula follows from the antisymmetric Cayley identity;
see \cite[Theorem~2.3 and Section~9]{CSS}.
For $N=2m+1$, L\H orincz, Raicu, Walther, and Weyman computed the
polynomial of the submaximal Pfaffian ideal
\cite[Theorem~3.9]{LRWW}. Their work also proves the strong
topological monodromy statement, including multiplicities, in that
case \cite[Section~5.2]{LRWW}. For general Pfaffian ideals,
Chen, Shi, Xu, and Zuo proved the monodromy and holomorphy conjectures
and described the resolution data \cite{CSXZ}. The first result below
determines the full polynomial in every nonmaximal case, and the
second strengthens the monodromy conclusion to the roots themselves
with the required multiplicity inequality.

For $r<m$, define the set of distinct values
\begin{align}
 R_{N,r}&=\left\{\frac{\binom{N-2i}{2}+u}{r-i}:
                 0\le i<r,\ 0\le u<r-i\right\},\label{eq:R}\\
 J_{N,r}(\rho)&=\left\{(i,u):0\le i<r,\ 0\le u<r-i,
              \ (r-i)\rho=\binom{N-2i}{2}+u\right\}.\label{eq:J}
\end{align}
The indices $i,u$ are integers. The distinction between the set of
values and its label fibers is essential.
\begin{maintheorem}[Roots and multiplicities]\label{thm:main}
If $1\le r<m$, the roots of $b_{\pr_r}(s)$ are exactly $-R_{N,r}$, and
\begin{equation}\label{eq:mainmult}
 \mult_{s=-\rho}b_{\pr_r}(s)=
 \begin{cases}
 |J_{N,r}(\rho)|,&\rho\in R_{N,r}\cap\Z,\\
 1,&\rho\in R_{N,r}\setminus\Z,\\
 0,&\rho\notin R_{N,r}.
 \end{cases}
\end{equation}
Equivalently,
\begin{equation}\label{eq:mainpoly}
 b_{\pr_r}(s)=
 \prod_{i=0}^{r-1}\left(s+\left\lceil\frac{\binom{N-2i}{2}}{r-i}\right\rceil\right)
 \prod_{\rho\in R_{N,r}\setminus\Z}(s+\rho),
\end{equation}
where the second product contains each distinct value once. If $r=m$, then
\begin{equation}\label{eq:boundary}
 b_{\pr_m}(s)=\prod_{j=0}^{m-1}(s+2j+1+2\epsilon).
\end{equation}
In particular $b_{\pr_1}(s)=s+D_0$ for every $N\ge2$.
\end{maintheorem}

The integer factors in \eqref{eq:mainpoly} have total degree $r$,
counted with multiplicity. Nonintegral collisions behave differently:
for $(N,r)=(21,6)$, each of $106/3$ and $107/3$ has two labels,
but each gives a simple root. Thus taking the product over all labels
would in general give an incorrect polynomial. Theorem~\ref{thm:main}
also determines every local polynomial by the transverse rank slice;
see Corollary~\ref{cor:local}.

Write $Z_{\pr_r}^{\mathrm{mot}}(T)$ and
$Z_{\pr_r}^{\mathrm{top}}(s)$ for the motivic and topological zeta
functions. The motivic series takes coefficients in
$\mathcal M_{\C}=K_0(\operatorname{Var}_{\C})[\mathbb L^{-1}]$,
where $\mathbb L=[\mathbb A^1]$. Its pole-order convention is specified
in Definition~\ref{def:motivicorder}; it does not treat
$\mathcal M_{\C}$ as a field.

\begin{maintheorem}[Strong motivic monodromy with multiplicities]
\label{thm:strong}
For every $1\le r\le m$, set
\[
 N_i=r-i,\qquad \nu_i=\binom{N-2i}{2}\quad(0\le i<r),
 \qquad p_{N,r}(\rho)=\#\{i:\nu_i/N_i=\rho\}.
\]
Then
\begin{equation}\label{eq:motivicdenominator}
 \left(\prod_{i=0}^{r-1}(1-\mathbb L^{-\nu_i}T^{N_i})\right)
 Z_{\pr_r}^{\mathrm{mot}}(T)\in\mathcal M_{\C}[T],
\end{equation}
and, for every positive rational number $\rho$,
\begin{equation}\label{eq:strongineq}
 \operatorname{ord}^{\mathrm{mot}}_{-\rho} Z_{\pr_r}^{\mathrm{mot}}
 \le p_{N,r}(\rho)
 \le \mult_{s=-\rho} b_{\pr_r}(s).
\end{equation}
In particular every motivic pole is a Bernstein--Sato root with at
least the same multiplicity. The analogous inequality holds for the
ordinary pole order of $Z_{\pr_r}^{\mathrm{top}}(s)$, and
\begin{equation}\label{eq:topdivisibility}
 b_{\pr_r}(s)Z_{\pr_r}^{\mathrm{top}}(s)\in\mathbb Q[s].
\end{equation}
\end{maintheorem}

The proof controls the common denominator before cancellation.
If two different indices have the same slope $\nu_i/N_i$, that
slope is an odd integer. The integral part of Theorem~\ref{thm:main}
therefore supplies every multiplicity needed in
\eqref{eq:strongineq}; a nonintegral slope occurs at most once.
This gives a uniform proof without requiring a noncancellation
argument for the motivic coefficients.

The multiplicity formula also has the following sharp consequence.
\begin{maintheorem}[Largest root multiplicity]\label{thm:maxmultiplicity}
For each $r\ge1$, the largest multiplicity of a root of $b_{\pr_r^{(N)}}$ as $N\ge2r+2$ varies is
\[
 M_r=\left\lfloor\frac{\sqrt{8r+1}-1}{2}\right\rfloor.
\]
The same maximum holds if the maximal boundaries $N=2r,2r+1$ are included, since their roots are simple. If $M_r\ge2$, it is attained for $N=2r+M_r^2$ at the root $-(4M_r^2-1)$.
\end{maintheorem}

The upper bound is obtained from finite relations in the graph module
and a skew Capelli recursion. A weighted Gaussian integral then forces
the total integral multiplicity to reach that upper bound. At a
nonintegral parameter, full central characters separate the layers of
a common-kernel filtration, so no nilpotent extension survives.
Root occurrence is proved separately using the homogeneous-map
criterion of \cite[Theorem~3.6]{CZanalytic} and transverse rank slices.
Appendix~\ref{sec:normal-cone-independent} gives a geometric explanation
of nonintegral simplicity: the nilpotent monodromy vanishes on the
whole nonunipotent Verdier specialization.

The paper is organized as follows. Section~\ref{sec:preliminaries}
fixes the graph-module convention, recalls highest weights and the
Pfaffian Schur decomposition, explains skew Capelli operators and
central-character tensor shifts, and defines the zeta functions and
their pole orders. Section~\ref{sec:upper} proves the balancing identity
for ordinary products, derives the bordering relations, and uses
Fourier transformation and the Capelli recursion to obtain the label
upper bound. Section~\ref{sec:analytic} computes the radial measure,
establishes the needed constant-term coefficient, constructs the
Mellin pole, and deduces the exact integral multiplicities.
Section~\ref{sec:central} identifies the central characters of the
common kernel and separates the filtration layers to prove
nonintegral simplicity. Section~\ref{sec:roots} proves the
homogeneous-map criterion, applies the rank slices, handles the maximal
cases, and completes the proof of Theorem~\ref{thm:main}.
Section~\ref{sec:multiplicity-monodromy} gives the local polynomials,
the interval formula for integral multiplicities, the sharp bound of
Theorem~\ref{thm:maxmultiplicity}, and the resolution proof of
Theorem~\ref{thm:strong}. In Appendix~\ref{sec:normal-cone-independent},
the exterior-power image, its finite cover, and the deformation to the
normal cone lead to an independent geometric proof of nonintegral
simplicity.

\begin{center}
\bf{Acknowledgements}
\end{center}

 AI tools (ChatGPT and DeepSeek) were used during the preparation of the paper as technical assistants to help with checking calculations and arguments, developing details of proofs, and improving the exposition. All mathematical statements, proofs, and conclusions were verified by the authors. Y. Chen and H. Zuo were supported by BJNSF Grant 1252009. H. Zuo was supported by NSFC Grant 12671056.

\section{Preliminaries}\label{sec:preliminaries}
\subsection{Pfaffians and the Bernstein--Sato convention}\label{sec:graph}
Let $S=\C[x_{ab}:1\le a<b\le N]$. Extend coordinates and derivatives
by $x_{ba}=-x_{ab}$, $\partial_{ba}=-\partial_{ab}$, and
$x_{aa}=\partial_{aa}=0$. For an ordered list
$I=(i_1,\ldots,i_{2t})$, define
\begin{equation}\label{eq:pfdef}
 \Pf_I(X)=\frac1{2^t t!}\sum_{\sigma\in\mathfrak S_{2t}}\sgn(\sigma)
 \prod_{j=1}^t x_{i_{\sigma(2j-1)},i_{\sigma(2j)}}.
\end{equation}
A repeated index gives zero, and $\Pf_\varnothing=1$.
We put $\pr_0=S$ and $\pr_t=0$ for $t>m$.

For a nonzero polynomial $f$, the classical polynomial $b_f(s)$ is
the monic polynomial of least degree satisfying
$P(s)f^{s+1}=b_f(s)f^s$ for some $P(s)\in\D_X[s]$.
For an ideal, we use the graph construction of
\cite[Sections~2.4--2.5]{BMS}, which is independent of the chosen
generators and specializes to this convention for a principal ideal.
Here $\D_X$ is the Weyl algebra generated by $S$ and its partial
derivatives, with $[\partial_{ab},x_{cd}]=\delta_{ac}\delta_{bd}$
for increasing index pairs and all other generator commutators zero.
Enumerate the generators of $\pr_r$ as $f_1,\ldots,f_a$, where $a=\binom N{2r}$. Add coordinates $t_1,\ldots,t_a$ and their derivatives. The graph module is the free $S$-module
\[
 M=S[\partial_{t_1},\ldots,\partial_{t_a}]\delta.
\]
For a multi-index $\alpha\in\N^a$, with $\N=\{0,1,2,\ldots\}$, its actions are
\begin{align}
 \partial_{x_\ell}(h\partial_t^\alpha\delta)
  &=(\partial_{x_\ell}h)\partial_t^\alpha\delta
   -\sum_{j=1}^a h(\partial_{x_\ell}f_j)\partial_t^{\alpha+e_j}\delta,\label{eq:grapha}\\
 t_j(h\partial_t^\alpha\delta)
  &=hf_j\partial_t^\alpha\delta-\alpha_jh\partial_t^{\alpha-e_j}\delta.\label{eq:graphb}
\end{align}
Terms with $\alpha_j=0$ in the second summand vanish. Give $t_j$ degree $1$, $\partial_{t_j}$ degree $-1$, and the $x$-Weyl operators degree zero. Let $V^\ell\D_{X,t}$ be spanned by normally ordered monomials of degree at least $\ell$. Put
\begin{equation}\label{eq:Qdef}
 U^\ell=V^\ell\D_{X,t}\delta,\qquad Q=U^0/U^1,\qquad s=-\sum_j\partial_{t_j}t_j.
\end{equation}
The graph definition of \cite{BMS} identifies $b_{\pr_r}$ with the minimal polynomial of $s$ on $Q$. The quotient is cyclic over $\operatorname{gr}_V^0\D_{X,t}$ and $s$ is central in this ring. Thus a polynomial in $s$ kills $Q$ if and only if it kills the class of $\delta$. Explicitly,
\begin{equation}\label{eq:saction}
 s(h\partial_t^\alpha\delta)=-\sum_j hf_j\partial_t^{\alpha+e_j}\delta+|\alpha|h\partial_t^\alpha\delta.
\end{equation}
For $E=\sum_{a<b}x_{ab}\partial_{ab}$, homogeneity gives $E\delta=rs\delta$.

The eigenvalues of $s$ on $Q$ are the roots of $b_{\pr_r}$; the
exponent of a factor $s+\rho$ is the nilpotence order of $s+\rho$
on its primary summand. We also use the canonical decreasing rational
$V$-filtration on $M$ along $t=0$. With this normalization,
$s+\beta$ is nilpotent on
$\operatorname{gr}_V^\beta M=V^\beta M/V^{>\beta}M$.
The filtration is discrete, exhaustive, and good over $V^0\D_{X,t}$;
$t_j V^\beta M\subseteq V^{\beta+1}M$ and
$\partial_{t_j}V^\beta M\subseteq V^{\beta-1}M$.
These conditions fix the index convention used in the appendix;
see \cite[Sections~1.1--1.3]{BMS}.

\subsection{Highest weights and the polynomial representation}
\label{subsec:representations}
Let $V=\C^N$, let $\mathfrak g=\mathfrak{gl}_N$, and let $e_{uv}$ denote its
matrix units. The action used in this section is the algebra homomorphism
\begin{equation}\label{eq:skew-g-action}
 \tau:U(\mathfrak g)\longrightarrow\D_X,
 \qquad \tau(e_{uv})=A_{uv}:=\sum_{a=1}^N x_{ua}\partial_{va}.
\end{equation}
Here both the coordinates and the derivatives are extended skew-symmetrically.
In particular,
\begin{equation}\label{eq:skew-action-commutators}
 [A_{uv},x_{ab}]=\delta_{va}x_{ub}+\delta_{vb}x_{au},
 \qquad
 \sum_{u=1}^N A_{uu}=2E.
\end{equation}
Thus the linear coordinate functions carry the representation $\bigwedge^2V$.
Write $Z=Z(U(\mathfrak g))$ for the center of the universal enveloping algebra;
explicitly,
\[
 U(\mathfrak g)=T(\mathfrak g)/\langle a\otimes b-b\otimes a-[a,b]:a,b\in\mathfrak g\rangle,
 \qquad T(\mathfrak g)=\bigoplus_{k\ge0}\mathfrak g^{\otimes k}.
\]
Its center consists of the elements commuting with all of
$U(\mathfrak g)$. A weight vector of weight $\lambda$ is a nonzero
simultaneous eigenvector of the diagonal matrix units, with
$e_{uu}w=\lambda_u w$. A highest-weight module is a module generated
by a highest-weight vector.

We use highest weights relative to the upper triangular matrices. A
highest-weight vector of weight $\lambda=(\lambda_1,\ldots,\lambda_N)$ is a
nonzero vector $w$ such that $e_{uv}w=0$ for $u<v$ and
$e_{uu}w=\lambda_u w$. The irreducible highest-weight module of weight $\lambda$
has a central character $\chi_\lambda:Z\to\C$: every $\zeta\in Z$
acts by the scalar $\chi_\lambda(\zeta)$. The Harish--Chandra theorem (see \cite[Section~1]{BG}) says that
\begin{equation}\label{eq:HC-multiset}
 \chi_\lambda=\chi_\mu
 \quad\Longleftrightarrow\quad
 \{\lambda_1+N-1,\lambda_2+N-2,\ldots,\lambda_N\}
 =\{\mu_1+N-1,\mu_2+N-2,\ldots,\mu_N\}
\end{equation}
as multisets. For fixed $\zeta$, its value $\chi_\lambda(\zeta)$ is a
polynomial in $\lambda$.

A partition $\lambda=(\lambda_1\ge\cdots\ge\lambda_m\ge0)$ has
size $|\lambda|=\sum_i\lambda_i$; write
$\lambda^{(2)}=(\lambda_1,\lambda_1,\ldots,\lambda_m,\lambda_m)$,
with a final zero when needed.
For a partition $\mu$ of length at most $N$, $\mathbb S_\mu\C^N$ denotes the irreducible
polynomial $\GL_N$-representation of highest weight $\mu$ (padded by zeros).
The variables transform as $\bigwedge^2\C^N$.
Finite-dimensional polynomial $\GL_N$-representations are completely reducible.
We use the classical multiplicity-free decomposition
\begin{equation}\label{eq:skew-polynomial-decomposition}
 S=\bigoplus_{\lambda}\mathbb S_{\lambda^{(2)}}\C^N,
 \qquad \deg \mathbb S_{\lambda^{(2)}}\C^N=|\lambda|,
\end{equation}
as recorded, for example, in \cite[(2.1)]{Perlman}.
If $c_j=\lambda'_j=\#\{i:\lambda_i\ge j\}$ are the column lengths,
a highest-weight vector of this constituent is
\[
 h_\lambda=\prod_{j=1}^{\lambda_1}\Pf_{[2c_j]}(X),
 \qquad [d]=\{1,\ldots,d\}.
\]
Indeed, each factor has weight $(1^{2c_j},0^{N-2c_j})$ and is killed
by upper triangular raising operators; their nonzero product has weight
$\lambda^{(2)}$.

\subsection{Skew Capelli operators and central characters}\label{subsec:capelli-prelim}
For $0\le p\le m$ define
\begin{equation}\label{eq:Cp}
 C_0=1,\qquad C_p=\sum_{|I|=2p}\Pf_I(X)\Pf_I(\partial),\qquad C_{m+1}=0.
\end{equation}
The index subsets are increasing and $C_1=E$.

We use the following established skew-Capelli facts. The image $\tau(Z)$ is
$\D_X^{\GL(V)}$, and this invariant algebra is generated by
$C_1,\ldots,C_m$. These statements are the Howe--Umeda theorem, restated
explicitly in \cite[p.~838]{Hashimoto}; see \cite{HU,KW}. The action in
\cite{Hashimoto} differs from ours by the automorphism
$e_{uv}\mapsto-e_{vu}$, which preserves the center and hence the assertion.

For $0\le i\le m$ and $a\in\C$, put
\[
 \lambda(i,a)=(\underbrace{a,\ldots,a}_{2i},
                  \underbrace{0,\ldots,0}_{N-2i}),
 \qquad
 c_p(i,a)=\binom{i}{p}\prod_{h=0}^{p-1}(a+2h),
 \quad 0\le p\le m,
\]
where $c_0(i,a)=1$ and $c_p(i,a)=0$ if $p>i$.

\begin{lemma}\label{lem:rectangle-character}
For every $\zeta\in Z$, choose a polynomial $P_\zeta$ such that
$\tau(\zeta)=P_\zeta(C_1,\ldots,C_m)$. Then
\begin{equation}\label{eq:rectangle-character}
 \chi_{\lambda(i,a)}(\zeta)
 =P_\zeta\bigl(c_1(i,a),\ldots,c_m(i,a)\bigr)
 \qquad(a\in\C).
\end{equation}
Consequently, on any subspace of a $\D$-module on which each $C_p$ acts by $c_p(i,a)$,
the whole center $Z$ acts by the exact character $\chi_{\lambda(i,a)}$.
\end{lemma}

\begin{proof}
For $a\in\N$, the polynomial $\Pf(X_{[2i]})^a$ is a highest-weight vector of
weight $\lambda(i,a)$. The rectangular skew-Capelli eigenvalue formula is
\[
 C_p\Pf(X_{[2i]})^a=c_p(i,a)\Pf(X_{[2i]})^a.
\]
It appears in \cite[Eq.~(3.10)]{KW}. It also follows directly from the
all-principal-minors Cayley identity \cite[Theorem~2.3]{CSS}: multiplication
by $\Pf(X_J)$ and summation over $|J|=2p$ give the factor $\binom{i}{p}$
by Pfaffian Laplace expansion. Terms differentiating a coordinate outside
$[2i]$ vanish.

Applying $\zeta$ to this vector proves \eqref{eq:rectangle-character} for
all nonnegative integral $a$. Both sides are polynomial in $a$, so the
identity holds for all $a\in\C$. This argument includes every element of
$\ker(\tau|_Z)$: its Harish--Chandra value on this line is identically zero.
It therefore proves the assertion for the full center, not merely for a
chosen collection of central lifts of the $C_p$. The last statement follows
by applying the polynomial identity $\tau(\zeta)=P_\zeta(C_1,\ldots,C_m)$
on the module in question.
\end{proof}

\subsection{Central support and tensor products}\label{subsec:central-support}
For a finite-dimensional algebraic $\GL_N$-module, its weights are
the integer tuples $\nu$ for which the corresponding simultaneous
diagonal eigenspace is nonzero. On a tensor product we use the
diagonal Lie algebra action
$x(\ell\otimes w)=x\ell\otimes w+\ell\otimes xw$.
The central support of a $Z$-module annihilated by a
finite-codimensional ideal is the finite set of characters appearing in its
generalized-character decomposition. Equivalently it is the zero set of
its annihilator in $\operatorname{Spec}(Z)$. We require the following form
of the central-character tensor theorem.

\begin{lemma}[Bernstein--Gelfand]\label{lem:central-tensor}
Let $M$ be a $U(\mathfrak g)$-module. Suppose an ideal of finite codimension in $Z$ annihilates $M$, and let $L$ be a finite-dimensional algebraic
$\GL_N$-representation. Then $L\otimes M$ is also annihilated by a
finite-codimensional ideal of $Z$. If the central support of $M$ is the
single character $\chi_\lambda$, then that of $L\otimes M$ is contained in
\[
 \{\chi_{\lambda+\nu}:\nu\text{ is a weight of }L\}.
\]
The assertion applies to arbitrary $M$; no finite-generation or
category-$\mathcal O$ hypothesis is required.
\end{lemma}

\begin{proof}
For a semisimple Lie algebra, \cite[Theorem~2.5]{BG} constructs the ideal of
relations between the two central actions on $L\otimes M$, universally for
all $M$. Its zero set is the union of the weight-shift correspondences.
By \cite[Corollary~2.6]{BG}, quotienting the input center by a
finite-codimensional annihilator of $M$ leaves a finite-dimensional algebra;
thus the output center has a finite-codimensional annihilator as well.
The support statement is stated immediately thereafter on p.~260.
For $\mathfrak{gl}_N=\mathfrak{sl}_N\oplus\C\operatorname{id}$, apply this
result to the semisimple part and retain the trace operator. On a trace
weight-$d$ summand of $L$, the trace operator on $L\otimes M$ is the trace
operator on $M$ plus $d$. Its generalized eigenvalue and nilpotence bound
therefore shift by exactly $d$. Combining this with the semisimple
weight shift gives $\lambda+\nu$ as asserted. The case $N=1$ is just this
trace calculation.
\end{proof}

\subsection{Motivic zeta functions and pole orders}\label{subsec:motivic}
The Grothendieck ring $K_0(\operatorname{Var}_{\C})$ is generated by
isomorphism classes $[Y]$ of complex algebraic varieties, with
$[Y]=[W]+[Y\setminus W]$ for a closed subvariety $W\subseteq Y$
and $[Y][Y']=[Y\times Y']$. Put
$\mathbb L=[\mathbb A^1]$ and
$\mathcal M_{\C}=K_0(\operatorname{Var}_{\C})[\mathbb L^{-1}]$.
For $X=\mathbb A^d$, an arc is a tuple in $\C[[t]]^d$; its
$q$-jet is its reduction modulo $t^{q+1}$. Let
\[
\mathcal X_{q,I}=\{\gamma\in J_q(X):
           \min_j\operatorname{ord}_t f_j(\gamma)=q\},
\]
where $J_q(X)=\mathbb A^{d(q+1)}$. We use the normalized,
non-equivariant motivic zeta function
\begin{equation}\label{eq:motivicdefinition}
Z_I^{\mathrm{mot}}(T)
=\sum_{q\ge0}[\mathcal X_{q,I}]\mathbb L^{-d(q+1)}T^q
\in\mathcal M_{\C}[[T]].
\end{equation}
The order of the zero truncated series is understood to be $>q$.
Multiplying this convention by $\mathbb L^d$, or omitting its
constant term, does not affect the denominator bounds below.

A log resolution of $I$ is a proper birational morphism
$\pi:Y\to X$ with $Y$ smooth such that
$I\mathcal O_Y=\mathcal O_Y(-\sum_{i\in A}N_iE_i)$ and the relevant
divisors have simple normal crossings: locally their union is a
union of coordinate hyperplanes. Write
$K_{Y/X}=\sum_{i\in A}(\nu_i-1)E_i$, and
\[
E_J^\circ=\left(\bigcap_{i\in J}E_i\right)
                \setminus\bigcup_{i\notin J}E_i,
\qquad E_\varnothing^\circ=Y\setminus\bigcup_{i\in A}E_i.
\]
For a principalization that is an isomorphism off $V(I)$, all the
relevant $N_i$ are positive. The change-of-variables formula for
motivic integration gives
\begin{equation}\label{eq:motivicresolution}
Z_I^{\mathrm{mot}}(T)=\mathbb L^{-d}
 \sum_{J\subseteq A}[E_J^\circ]
 \prod_{i\in J}
 \frac{(\mathbb L-1)\mathbb L^{-\nu_i}T^{N_i}}
      {1-\mathbb L^{-\nu_i}T^{N_i}}.
\end{equation}
We use this standard resolution formula in the ideal version
recorded in \cite[Section~2.5, (2.1)]{CSXZ}; for the underlying motivic
change of variables see \cite[Lemma~3.3]{DL}. The factor $\mathbb L^{-d}$
corresponds to the normalization in \eqref{eq:motivicdefinition}.
The topological zeta function is
\begin{equation}\label{eq:topologicaldefinition}
Z_I^{\mathrm{top}}(s)=
\sum_{J\subseteq A}\chi_c(E_J^\circ)
                    \prod_{i\in J}\frac1{N_i s+\nu_i},
\end{equation}
where $\chi_c$ is the compactly supported Euler characteristic.
Both formulas are independent of the chosen resolution.
The parameter $s$ is related formally to $T$ by $T=\mathbb L^{-s}$,
so the factor $1-\mathbb L^{-a}T^b$ has candidate slope $s=-a/b$.

\begin{definition}[Denominator convention for motivic pole orders]
\label{def:motivicorder}
Consider a rational motivic series expressible as a finite sum
\begin{equation}\label{eq:denominatorpresentation}
Z(T)=\sum_{\alpha}P_\alpha(T)
          \prod_{j=1}^{h_\alpha}
           (1-\mathbb L^{-a_{\alpha j}}T^{b_{\alpha j}})^{-1},
\qquad
P_\alpha(T)\in\mathcal M_{\C}[T],\quad
 a_{\alpha j},b_{\alpha j}\in\Z_{>0}.
\end{equation}
For a positive rational number $\rho$, its pole order at
$s=-\rho$, denoted by $\operatorname{ord}^{\mathrm{mot}}_{-\rho}Z$, is the least nonnegative integer $e$ for which such a
presentation exists with at most $e$ factors satisfying
$a_{\alpha j}/b_{\alpha j}=\rho$ in each summand, repetitions counted.
A positive order means that $-\rho$ is a pole. This is a denominator
notion of order in the motivic ring. After Hodge realization, with $\mathbb L$ sent to the indeterminate
$uv$, the ordinary order along $T=(uv)^\rho$ (in a suitable finite
field extension) is bounded by this order.
\end{definition}

This convention states an upper bound without assuming unique
factorization in $\mathcal M_{\C}$. Formula~\eqref{eq:motivicresolution}
shows that the order at $-\rho$ is at most the number of divisors
with $\nu_i/N_i=\rho$. Theorem~\ref{thm:strong} will prove the
stronger assertion that one common denominator works simultaneously
with the required bounds at every slope. For
\eqref{eq:topologicaldefinition}, pole orders have their usual meaning
for a rational function of the complex variable $s$.

\section{An upper bound for the roots and their multiplicities}\label{sec:upper}
We construct an annihilating polynomial for the graph quotient. This
both excludes roots outside the stated windows and bounds the
multiplicity of each root by the number of its labels.

\subsection{Balancing ordinary products}\label{subsec:balancing}
\begin{lemma}[Balancing two Pfaffian factors]\label{lem:balancing}
For integers $0\le a<b\le m$ with $b\ge a+2$,
\begin{equation}\label{eq:balancing}
 \pr_a\pr_b\subseteq\pr_{a+1}\pr_{b-1}.
\end{equation}
More generally, if $0\le c_1,\ldots,c_v\le m$ and
$\sum_jc_j=rv$ for an integer $r$, then
\begin{equation}\label{eq:balanced-product}
 \prod_{j=1}^v\pr_{c_j}\subseteq\pr_r^v.
\end{equation}
\end{lemma}
\begin{proof}
The degree-$(a+b)$ generators of the first product are the image of
$\bigwedge^{2a}\C^N\otimes\bigwedge^{2b}\C^N$ under actual polynomial
multiplication. The two-column Pieri rule gives the constituents
$\mathbb S_{(2^h,1^{2a+2b-2h})}\C^N$, $0\le h\le2a$, omitting
partitions of length greater than $N$. Write $e_j$ and $s_\mu$ for the
elementary symmetric and Schur characters. One can read this rule directly
from the dual Jacobi--Trudi identity
\[
 s_{(2^h,1^{p+q-2h})}=e_{p+q-h}e_h-e_{p+q-h+1}e_{h-1}:
\]
summing over $0\le h\le p\le q$, with $e_{-1}=0$, telescopes to
$e_pe_q$, the character of the two exterior powers.
In \eqref{eq:skew-polynomial-decomposition}, paired rows force $h=2j$.
Thus every constituent that can occur in the multiplication image has
\[
 \lambda=(2^j,1^{a+b-2j}),\qquad
 0\le j\le a,\qquad a+b-j\le m,
\]
and has highest-weight vector
$h_j=\Pf_{[2j]}\Pf_{[2(a+b-j)]}$.
We will put every such vector in the target ideal; no assertion that
every tensor constituent survives multiplication is being assumed.

Fix a feasible $j$. Set $A=[2j]$, $B=[2(a+b-j)]$, $R=B\setminus A$,
$t=a+b-2j$ and $k=a+1-j$. Both $k$ and $t-k=b-1-j$ are positive.
For $K\subset R$ of size $2k$, let $\sigma(K,R\setminus K)$ be the
sign of the permutation listing first the increasing elements of $K$
and then the increasing elements of $R\setminus K$, relative to the
increasing order on $R$. We claim
\begin{equation}\label{eq:pfaffian-balancing-identity}
 \sum_{\substack{K\subset R\\|K|=2k}}
 \sigma(K,R\setminus K)\Pf_{A\cup K}\Pf_{A\cup(R\setminus K)}
 =\binom tk\Pf_A\Pf_B.
\end{equation}
For $j>0$, work first on the dense open set $\Pf_A\ne0$.
Writing the skew matrix on $B$ in blocks
$\left(\begin{smallmatrix}X_A&U\\-U^t&D\end{smallmatrix}\right)$,
put $T=D+U^tX_A^{-1}U$.
The skew Schur-complement identity gives
$\Pf_{A\cup K}=\Pf_A\Pf(T_K)$ and
$\Pf_B=\Pf_A\Pf(T)$.
For the two-form $\omega_T$ corresponding to $T$, the identity
\[
 \frac{\omega_T^k}{k!}\wedge\frac{\omega_T^{t-k}}{(t-k)!}
 =\binom tk\frac{\omega_T^t}{t!}
\]
is exactly the required sum of products of Pfaffians of $T$.
Multiplying by $\Pf_A^2$ proves
\eqref{eq:pfaffian-balancing-identity} on that dense open set, hence
as a polynomial identity everywhere. For $j=0$, the same exterior
identity proves it directly, with $\Pf_A=1$.
The factors on its left have half-sizes $a+1$ and $b-1$.
Since $\binom tk\ne0$ over $\C$, every $h_j$ lies in
$\pr_{a+1}\pr_{b-1}$.
This target is $\GL_N$-stable, and the ambient polynomial ring is
multiplicity free. It therefore contains every constituent that the
source multiplication image can have. This proves inclusion of their
degree-$(a+b)$ generating spaces and hence \eqref{eq:balancing} for
the ordinary ideals.

For \eqref{eq:balanced-product}, whenever two half-sizes have difference
at least two, apply \eqref{eq:balancing} to a smallest and a largest
one. Their sum is unchanged and their sum of squares decreases by
$2(b-a-1)>0$. The process terminates with all half-sizes differing by
at most one. Their average is the integer $r$, so they all equal $r$.
Multiplying each inclusion by the other ideal factors proves the result.
\end{proof}

\begin{corollary}[The lowest degree of an ordinary power]\label{cor:ordinary-power}
For every $v\ge0$,
\begin{equation}\label{eq:powers}
 (\pr_r^v)_{rv}=
 \bigoplus_{\substack{|\lambda|=rv\\\lambda_1\le v}}
 \mathbb S_{\lambda^{(2)}}\C^N.
\end{equation}
\end{corollary}
\begin{proof}
For $v=0$ both sides are the constants. For $v>0$, the generators of
the power are products of $v$ vectors in $\bigwedge^{2r}\C^N$.
Every coordinate of a weight of one exterior power is zero or one;
therefore every weight coordinate in this tensor product is at most
$v$. Decomposition \eqref{eq:skew-polynomial-decomposition} proves
the inclusion into the displayed right side.
Conversely, take $|\lambda|=rv$, $\lambda_1\le v$, and pad its column
lengths $c_j$ by zeros to a list of length $v$.
Their sum is $rv$, and $h_\lambda\in\prod_j\pr_{c_j}$.
Lemma~\ref{lem:balancing} puts this nonzero highest-weight vector in
$\pr_r^v$. Stability under $\GL_N$ gives the entire constituent.
This proves the reverse inclusion.
\end{proof}
Formula \eqref{eq:powers} is the ordinary-power formula of
Abeasis--Del Fra \cite{ADF}, recorded in \cite[(2.6)]{Perlman}.
The proof also gives the product inclusion needed for the graph relations.

\subsection{Bordering operators and graph relations}\label{subsec:bordering}
For $a<b$ set
\begin{equation}\label{eq:Rop}
 R_{ab}=\sum_{i<j}(x_{ia}x_{jb}-x_{ib}x_{ja})\partial_{ij},\qquad
 F_{ab}=R_{ab}-sx_{ab},\qquad b_{ab}\Pf_I=\Pf_{I,a,b}.
\end{equation}
Bordering acts on ordered Pfaffians; two pair-borderings commute because exchanging two blocks of length two has positive sign.

\begin{lemma}\label{lem:border}
For every ordered Pfaffian,
\begin{equation}\label{eq:borderid}
 R_{ab}\Pf_I=x_{ab}\Pf_I-b_{ab}\Pf_I.
\end{equation}
The $F_{ab}$ commute. If $\eta_1,\ldots,\eta_k$ are entry pairs, then in $M$
\begin{equation}\label{eq:partition}
 F_{\eta_1}\cdots F_{\eta_k}\delta
 =\sum_\pi\sum_{(j_B)_{B\in\pi}}(-1)^{k-|\pi|}
 \left(\prod_{B\in\pi}b_Bf_{j_B}\right)
 \left(\prod_{B\in\pi}\partial_{t_{j_B}}\right)\delta.
\end{equation}
Here $\pi$ runs through the set partitions of $\{1,\ldots,k\}$, $1\le j_B\le a$, and $b_B=\prod_{\ell\in B}b_{\eta_\ell}$. The empty partition gives $\delta$ when $k=0$.
\end{lemma}
\begin{proof}
Expanding~\eqref{eq:pfdef} by the pair of additional indices proves~\eqref{eq:borderid}: the term pairing $a$ with $b$ is $x_{ab}\Pf_I$, and the other terms are the negative of the differential of $\Pf_I$ in the direction with entries $x_{ia}x_{jb}-x_{ib}x_{ja}$. This proof also covers repeated indices.

For fixed skew matrices $A,B$, let $R_A X=-XAX$. Differentiating $R_BX=-XBX$ in the direction $R_AX$ gives
\[
 XAXBX+XBXAX,
\]
which is symmetric in $A,B$. Hence $[R_A,R_B]=0$ in every size. Taking the skew matrix units gives the claimed commutation of the $R_{ab}$. Moreover $R_{ab}x_{cd}=R_{cd}x_{ab}$. Since $s$ commutes with the $x$-Weyl algebra, this proves $[F_{ab},F_{cd}]=0$.

Combining~\eqref{eq:borderid} with the graph actions gives
\[
 F_\eta(h\partial_t^\alpha\delta)
 =(R_\eta h-|\alpha|x_\eta h)\partial_t^\alpha\delta
 +\sum_j h(b_\eta f_j)\partial_t^{\alpha+e_j}\delta.
\]
In a term of~\eqref{eq:partition} with $v$ blocks, its coefficient is a product of $v$ Pfaffians. Applying $R_\eta$ gives $vx_\eta$ times that product minus the terms that add the new border to one existing block. The first contribution cancels $-|\alpha|x_\eta$ since $|\alpha|=v$. The graph-derivative term creates a singleton block. These are exactly the set-partition recursion and its signs, proving the formula by induction.
\end{proof}

\begin{proposition}\label{prop:graphrels}
For $0\le k\le r$ and all entry pairs $\eta_1,\ldots,\eta_k$, one has in $Q$
\begin{equation}\label{eq:graphrels}
 \pr_{r-k}(X)F_{\eta_1}\cdots F_{\eta_k}\delta=0.
\end{equation}
In particular every length-$r$ word in the $F_{ab}$ kills $\delta$.
\end{proposition}
\begin{proof}
For $k=0$, $f_j\delta=t_j\delta\in U^1$. For $k>0$, consider a term of~\eqref{eq:partition} with $v$ blocks of sizes $j_1,\ldots,j_v$, where $\sum j_\ell=k$. Its coefficient lies in $\pr_{r+j_1}\cdots\pr_{r+j_v}$. After multiplication by an $(r-k)$-Pfaffian, the $v+1$ half-sizes sum to
$(r-k)+\sum_{\ell=1}^v(r+j_\ell)=r(v+1)$. If any size exceeds $m$ the term is zero; otherwise the balancing lemma gives
\[
 \pr_{r-k}\pr_{r+j_1}\cdots\pr_{r+j_v}\subseteq\pr_r^{v+1}.
\]
For $k=r$ the first factor is $\pr_0=S$, explicitly allowed in the balancing lemma. Express such a coefficient as a sum of polynomial multiples of $f^\beta$ with $|\beta|=v+1$. Since the auxiliary-derivative degree is $|\alpha|=v$,
\[
 f^\beta\partial_t^\alpha\delta=\partial_t^\alpha t^\beta\delta\in U^1.
\]
All the terms therefore vanish in $Q$.
\end{proof}

\subsection{Fourier transformation}\label{subsec:fourier}
Apply the Weyl algebra automorphism $x_{ab}\mapsto\partial_{ab}$, $\partial_{ab}\mapsto-x_{ab}$, $s\mapsto s$ to the $\D_X[s]$-relations of $\delta$ in $Q$, and denote the resulting cyclic vector by $v$. Put
\begin{equation}\label{eq:qz}
 q=-rs-D_0,\qquad z=s+2N-2,\qquad
 B_{ab}=\sum_{i<j}x_{ij}(\partial_{ia}\partial_{jb}-\partial_{ib}\partial_{ja})+z\partial_{ab}.
\end{equation}
Normal ordering gives the Fourier image $F_{ab}\mapsto-B_{ab}$: the correction in
\[
 \sum_{i<j}(\partial_{ia}\partial_{jb}-\partial_{ib}\partial_{ja})x_{ij}
 =\sum_{i<j}x_{ij}(\partial_{ia}\partial_{jb}-\partial_{ib}\partial_{ja})+2(N-1)\partial_{ab}
\]
comes from the two derivative placements and the $N-1$ indices meeting a fixed index. It may also be checked directly using $[\partial_{ab},x_{ij}]=\delta_{ai}\delta_{bj}-\delta_{aj}\delta_{bi}$ with skew extension. Consequently the $B_{ab}$ commute and $[E,B_{ab}]=-B_{ab}$. Proposition~\ref{prop:graphrels} becomes
\begin{equation}\label{eq:fourierrels}
 Ev=qv,\quad EB^kv=(q-k)B^kv,\quad \pr_{r-k}(\partial)B^kv=0\ (0\le k\le r),\quad B^rv=0.
\end{equation}
Here $B^k$ denotes any word of length $k$, not a single selected operator.

\subsection{A skew Capelli recursion}\label{subsec:capelli-recursion}
\begin{lemma}[Skew Capelli recursion]\label{lem:capelli}
For $0\le p\le m$,
\begin{equation}\label{eq:capelli}
 \sum_{a<b}[C_p,x_{ab}]B_{ab}
 =C_p\bigl(E+pz-p(2p-1)\bigr)-(p+1)C_{p+1}.
\end{equation}
In particular, if every $B_{ab}$ kills $w$ in a $\D_X[s]$-module, then
\begin{equation}\label{eq:capelliproduct}
 t!C_tw=\prod_{i=0}^{t-1}\bigl(E+iz-i(2i-1)\bigr)w\qquad(0\le t\le m).
\end{equation}
\end{lemma}
\begin{proof}
Use independent commuting skew variables $\Xi=(\xi_{ab})$ for normal symbols, with coefficients in $X$ placed left of derivatives. Write
\[
 P_p=\sum_{|I|=2p}\Pf_I(X)\Pf_I(\Xi),\quad e=\sum_{a<b}x_{ab}\xi_{ab},\quad
 K_{ab}=\sum_{i<j}x_{ij}(\xi_{ia}\xi_{jb}-\xi_{ib}\xi_{ja}).
\]
As matrices $K=-\Xi X\Xi$. The Pfaffian minor-summation identity gives, for $A=X\Xi^t=-X\Xi$,
\begin{equation}\label{eq:gensquare}
 \det(1+tA)=\left(\sum_{p=0}^m P_pt^p\right)^2.
\end{equation}
One direct proof expands the Pfaffian of $\left(\begin{smallmatrix}tX&1\\-1&-\Xi\end{smallmatrix}\right)$: matching the remaining indices across the identity blocks leaves equal subsets in $X$ and $\Xi$, yielding $(-1)^{N(N-1)/2}\sum P_pt^p$. Squaring this Pfaffian and taking the block determinant gives~\eqref{eq:gensquare}.

Under $\dot\Xi=K$ one has $\dot A=A^2$. Differentiating~\eqref{eq:gensquare}, or comparing the elementary symmetric Newton identities for the paired eigenvalues of $A$, gives
\begin{equation}\label{eq:firstsymbol}
 \sum_{a<b}(\partial_{\xi_{ab}}P_p)K_{ab}=eP_p-(p+1)P_{p+1}.
\end{equation}
For an entirely formal derivation, put $P(t)=\sum P_pt^p$. Since $P(0)=1$, equation~\eqref{eq:gensquare} gives
\[
 \log P(t)=\frac12\tr\log(1+tA).
\]
Using $\dot A=A^2$ in this power-series identity yields $\dot P(t)=eP(t)-P'(t)$. Comparing coefficients proves~\eqref{eq:firstsymbol}, without any diagonalizability assumption.

The normal-symbol product rule is $\sigma(UV)=\sum_\alpha(\partial_\Xi^\alpha\sigma(U))(\partial_X^\alpha\sigma(V))/\alpha!$. Since $K$ is linear in $X$, the only extra term here is
\begin{equation}\label{eq:secondsymbol}
 \sum_{a<b,\,i<j}(\partial_{\xi_{ij}}\partial_{\xi_{ab}}P_p)
 (\xi_{ia}\xi_{jb}-\xi_{ib}\xi_{ja})=-2p(p-1)P_p.
\end{equation}
Here is a verification of the correction that specifies the signs. Let $L$ be the differential operator in $\Xi$ given by the left contraction in~\eqref{eq:secondsymbol}. Its normal symbol in $\Xi$ and a dual skew variable $H=(\eta_{ab})$ is
\[
 \left(\sum_{a<b}\xi_{ab}\eta_{ab}\right)^2
 -2\sum_{|I|=4}\Pf_I(\Xi)\Pf_I(H),
\]
by~\eqref{eq:firstsymbol} for $p=1$ with the variable roles exchanged. Both terms are invariant contractions. Normal quantization commutes with linear changes of variables, so $L$ is $\GL_N$-equivariant. The multiplicity-free decomposition used in Section~\ref{subsec:representations} shows that it acts by a scalar on the span of the $2p$-Pfaffians. Evaluate this scalar on one Pfaffian at $\Xi=\operatorname{diag}(J,\ldots,J)$, where $J=\left(\begin{smallmatrix}0&1\\-1&0\end{smallmatrix}\right)$ and the Pfaffian is one. Each unordered pair of canonical blocks contributes four cross-edge terms, all $-1$: for the first two blocks this follows directly from
\[
 \Pf_{1234}(\Xi)=\xi_{12}\xi_{34}-\xi_{13}\xi_{24}+\xi_{14}\xi_{23},
\]
using the ordered derivative pairs $(13,24),(24,13),(14,23),(23,14)$ and their coefficients in $L$. The other blocks contribute one. Same-edge second derivatives vanish, and the remaining terms vanish at this point. Thus the scalar is $-4\binom p2=-2p(p-1)$, proving~\eqref{eq:secondsymbol}.

The $z$ term contributes $pzP_p$ by homogeneity, whereas the normal symbol of $C_pE$ is $eP_p+pP_p$. Combining these facts with~\eqref{eq:firstsymbol}--\eqref{eq:secondsymbol} proves~\eqref{eq:capelli}, including $p=0$. If $Bw=0$, its left side vanishes. Each $C_p$ commutes with $E$ and $z$; induction starting with $C_0=1$ gives~\eqref{eq:capelliproduct}.
\end{proof}

\subsection{The complete label bound}\label{subsec:label-bound}
\begin{proposition}\label{prop:upper}
For every $1\le r\le m$,
\begin{equation}\label{eq:upper}
 b_{\pr_r}(s)\ \mid\ 
 \prod_{i=0}^{r-1}\prod_{u=0}^{r-i-1}
 \bigl((r-i)s+\tbinom{N-2i}{2}+u\bigr).
\end{equation}
\end{proposition}
\begin{proof}
Set $P_k(q,z)=\prod_{i=0}^{r-k-1}(q-k+iz-i(2i-1))$. Start with $B^rv=0$ and descend through $k=r-1,\ldots,0$. Suppose $w$ is a polynomial in $s$ times $v$ and every length-$(k+1)$ word kills $w$. If $u=B^kw$, then
\[
 B_{ab}u=0,\qquad C_{r-k}u=0,\qquad Eu=(q-k)u,
\]
by~\eqref{eq:fourierrels}. Lemma~\ref{lem:capelli} gives $P_k(q,z)u=0$. Replacing $w$ by $P_k(q,z)w$ makes all length-$k$ words vanish. At the end, $\prod_{k=0}^{r-1}P_k(q,z)$ kills $v$. The factors are
\[
 q-k+iz-i(2i-1)=-\bigl((r-i)s+\tbinom{N-2i}{2}+k\bigr).
\]
The index sets $0\le k<r$, $0\le i<r-k$ and $0\le i<r$, $0\le k<r-i$ coincide. Inverse Fourier transformation kills $\delta$ in $Q$, and centrality plus cyclicity kills the whole quotient. Its minimal polynomial therefore divides the displayed product.
\end{proof}

\begin{remark}\label{rem:maxupper}
For $r=m$, all bordered generators vanish identically, so already $F_{ab}\delta=0$ and $B_{ab}v=0$. Applying~\eqref{eq:capelliproduct} with $C_mv=0$ gives the sharper upper bound with $u=0$ only. Since $N-2i=2(m-i)+\epsilon$,
\[
 \frac{\binom{N-2i}{2}}{m-i}=2(m-i)-1+2\epsilon.
\]
Thus the upper bound is the simple product in~\eqref{eq:boundary}. Proposition~\ref{prop:even-boundary} proves attainment in even size; the odd-size result is \cite[Theorem 3.9]{LRWW}.
\end{remark}

\section{Integral multiplicities from a weighted Gaussian integral}
\label{sec:analytic}

We construct a Mellin integral with a pole of order $r$ at an
integer parameter. A Bernstein equation converts this pole into a lower
bound for the total integral multiplicity. Combined with
Proposition~\ref{prop:upper}, it determines each integral multiplicity.

\subsection{Radial coordinates and the Jacobian}\label{subsec:radial}

Write $N=2m+\epsilon$, where $\epsilon\in\{0,1\}$, and put
$D_0=\binom N2$.  We identify the space of skew-symmetric matrices with
$\C^{D_0}$ using the coordinates $x_{ij}$, $i<j$.  Its Lebesgue measure is
normalized by $\int_{\C}e^{-|z|^2}\,dz=\pi$.  For $1\leq j\leq m$ set
\[
 Q_j(X)=\sum_{\substack{A\subseteq\{1,\ldots,N\}\\|A|=2j}}
             |\Pf(X_A)|^2,
 \qquad \|X\|^2=\sum_{i<j}|x_{ij}|^2.
\]
Thus $\pr_j$ denotes the ideal generated by the Pfaffians appearing in $Q_j$.
For real variables $\lambda_1,\ldots,\lambda_m$, let $e_j(\lambda)$ be the
$j$th elementary symmetric polynomial and let
$\Delta_m(\lambda)=\prod_{i<j}(\lambda_i-\lambda_j)$.
A Schwartz function is a smooth function all of whose derivatives decrease
faster than every inverse polynomial.

\begin{lemma}[Radial measure]\label{lem:skew-radial}
Put
\[
 J=\begin{pmatrix}0&1\\-1&0\end{pmatrix},\qquad
 X_\lambda=\operatorname{diag}(\sqrt{\lambda_1}J,\ldots,\sqrt{\lambda_m}J,0_\epsilon).
\]
There is a constant $c_N>0$ such that every integrable
unitary-congruence-invariant function $F$ satisfies the radial formula
\[
 \int_{\C^{D_0}}F(X)\,dX
 =c_N\int_{\lambda_1>\cdots>\lambda_m>0}
 F(X_\lambda)\Delta_m(\lambda)^4
 \prod_{i=1}^m\lambda_i^{2\epsilon}\,d\lambda.
\]
On this diagonal form, $Q_j=e_j(\lambda)$ and
$\|X\|^2=\sum_i\lambda_i$.
\end{lemma}

\begin{proof}
We first recall the normal form and give the needed Jacobian computation.
The antilinear map $T(v)=X\overline v$ satisfies $T^2=-XX^*$.
For the Hermitian inner product linear in the second argument,
skew-symmetry gives
$\langle u,Tv\rangle=-\langle v,Tu\rangle$.
Choose a unit eigenvector $u$ of $XX^*$ with eigenvalue $\sigma^2>0$.
Then $\|Tu\|=\sigma$ and $\langle u,Tu\rangle=0$.
Putting $v=-Tu/\sigma$ gives an orthonormal pair with
$Tu=-\sigma v$ and $Tv=\sigma u$.
Its orthogonal complement is $T$-stable by the same inner-product identity.
Induction produces the asserted normal form under $X\mapsto UXU^{\mathsf T}$,
with $U$ unitary.  On the full-rank locus the remaining zero block has size
$\epsilon$.

At $D=\operatorname{diag}(\sigma_1J,\ldots,\sigma_mJ,0_\epsilon)$,
an angular tangent vector is
$AD+DA^{\mathsf T}$, where $A^*=-A$.
For $i<j$, its $(i,j)$ block, with $B=A_{ij}\in M_2(\C)$, is
\[
 C=\sigma_jBJ-\sigma_iJ\overline B.
\]
Right multiplication by $-J$ turns this into
$\sigma_jB+\sigma_iS(B)$, where $S(B)=J\overline B J$.
The real-linear map $S$ is an involution.  Its positive eigenspace consists
of the matrices
$\bigl(\begin{smallmatrix}a&b\\\overline b&-\overline a\end{smallmatrix}\bigr)$,
and its negative eigenspace of
$\bigl(\begin{smallmatrix}a&b\\-\overline b&\overline a\end{smallmatrix}\bigr)$.
Each has real dimension four.  The corresponding absolute real determinant
is therefore
\[
 |\sigma_j+\sigma_i|^4|\sigma_j-\sigma_i|^4
       =|\lambda_j-\lambda_i|^4.
\]
If $\epsilon=1$, the block joining the $i$th two-plane to the last
coordinate is a real-linear map on $\C^2$ with absolute determinant
$\sigma_i^4=\lambda_i^2$.
The radial and phase directions of each diagonal block contribute a
constant times $\sigma_i\,d\sigma_i=\tfrac12d\lambda_i$.
Integrating the remaining angular variables over their compact homogeneous
space gives a positive constant independent of $\lambda$.
The omitted lower-rank and repeated-singular-value loci have measure zero.
This proves the measure formula; its homogeneity also checks since
$m+2\epsilon m+2m(m-1)=D_0$.

Finally, the coordinates of
$\bigl(\sum_{i<j}x_{ij}e_i\wedge e_j\bigr)^j/j!$
in the orthonormal exterior basis are the principal Pfaffians $\Pf(X_A)$.
The induced action of a unitary matrix on the exterior power is unitary.
Evaluation at the block diagonal matrix gives
$Q_j=e_j(\sigma_1^2,\ldots,\sigma_m^2)$, as claimed.
\end{proof}

\subsection{A constant-term coefficient}\label{subsec:constant-term}
We will need a particular coefficient of a fourth Vandermonde power.
For a Laurent polynomial, $\operatorname{CT}$ denotes its constant term;
for a power series, $[z^\alpha]$ denotes coefficient extraction.

\begin{theorem}[Morris constant-term identity]\label{thm:morris}
For nonnegative integers $a,b,k$ and a positive integer $p$,
\[
 \operatorname{CT}_{z_1,\ldots,z_p}
 \prod_{i=1}^p(1-z_i)^a(1-z_i^{-1})^b
 \prod_{i\ne j}(1-z_i/z_j)^k
 =\prod_{j=0}^{p-1}
 \frac{(a+b+kj)!\,(k(j+1))!}
 {(a+kj)!\,(b+kj)!\,k!}.
\]
\end{theorem}
This is the classical identity, in the normalization of
\cite[Theorem~1.2]{Morris}.

\begin{lemma}\label{lem:fourth-vandermonde}
For integers $p\geq1$ and $h\geq0$, put $H=h+2(p-1)$.  Then
\[
 [z_1^H\cdots z_p^H]\Delta_p(z)^4
                (z_1+\cdots+z_p)^{ph}
 =(ph)!\prod_{j=0}^{p-1}\frac{(2j+2)!}{2!(h+2j)!}>0.
\]
\end{lemma}

\begin{proof}
In Theorem~\ref{thm:morris}, set $b=h$ and $k=2$.
Divide the identity by $a^{ph}$ and let the positive integer $a$ tend to
infinity.  A constant term is unchanged by replacing every $z_i$ by $z_i/a$.
In the resulting coefficient extraction only finitely many coefficients
of the positive-power factors are relevant, with a bound independent of
$a$.  Consequently its limit is
\[
 \operatorname{CT}\prod_{i=1}^p e^{-z_i}(-z_i^{-1})^h
                  \prod_{i\ne j}(1-z_i/z_j)^2.
\]
Now
$\prod_{i\ne j}(1-z_i/z_j)^2
 =\Delta_p(z)^4/\prod_i z_i^{2(p-1)}$.
By homogeneity only the degree-$ph$ part of $e^{-\sum z_i}$ contributes.
Its sign cancels the sign $(-1)^{ph}$ from the inverse powers.
Thus the displayed limit is the coefficient in the lemma divided by
$(ph)!$.
On the product side, each ratio
$(a+h+2j)!/(a+2j)!$ is a monic degree-$h$ polynomial in $a$.
Taking leading coefficients proves the formula and its nonvanishing.
\end{proof}

\subsection{A Mellin pole of order \texorpdfstring{$r$}{r}}\label{subsec:mellin}
\begin{proposition}\label{prop:pfaffian-mellin}
Suppose $2\leq r<m$.  Define
\[
 \phi_r(X)=Q_m(X)^{2r}e^{-\|X\|^2},
 \qquad
 \rho=(m-r+1)(2m+2\epsilon+1).
\]
The meromorphic continuation of
\[
 Z_r(s)=\int_{\C^{D_0}}Q_r(X)^s\phi_r(X)\,dX
 \qquad (\operatorname{Re}s>0)
\]
has a pole of exactly order $r$ at $s=-\rho$.
\end{proposition}

\begin{proof}
The weight $\phi_r$ is a polynomial in the real coordinates times a
Gaussian and is therefore Schwartz.  Lemma~\ref{lem:skew-radial} gives
\[
 Z_r(s)=c_N\int_{\lambda_1>\cdots>\lambda_m>0}
 e_r(\lambda)^s\Delta_m(\lambda)^4
       \prod_i\lambda_i^{2\epsilon+2r}e^{-\sum_i\lambda_i}\,d\lambda.
\]
Use ratio coordinates $t_1=\lambda_1$ and
$t_j=\lambda_j/\lambda_{j-1}$ for $j\geq2$, so
$t_1>0$ and $0<t_j<1$ for $j\geq2$.
Set
\begin{align*}
 L&=1+t_2+t_2t_3+\cdots+t_2\cdots t_m,\\
 V(t)&=\prod_{1\leq i<j\leq m}(1-t_{i+1}\cdots t_j)^4,
 &U(t)&=\frac{e_r(\lambda)}{\lambda_1\cdots\lambda_r},\\
 a_k(s)&=(m-k)\bigl(2(m-k)-1+2\epsilon+2r\bigr)
                   +(r-k)_+s\quad (0\leq k<m).
\end{align*}
The quotient defining $U$ is a polynomial in $t_2,\ldots,t_m$ with
nonnegative coefficients and constant term one: a term selected at
$i_1<\cdots<i_r$ has $i_j\geq j$.
Hence $1\leq U\leq\binom mr$ and $1\leq L\leq m$ on the closed cube.
Combining the ratio-coordinate Jacobian with the displayed density,
the exponent of $t_j$ is $a_{j-1}(s)-1$.
Integration of $t_1$ therefore yields
\begin{equation}\label{eq:pfaffian-mellin-coordinates}
 Z_r(s)=c_N\Gamma(a_0(s))
 \int_{(0,1)^{m-1}}\prod_{j=2}^m t_j^{a_{j-1}(s)-1}
              L^{-a_0(s)}V(t)U(t)^s\,dt.
\end{equation}

Put $p=r-1$ and $h=(m-r)(2m+2\epsilon+1)-2m$.
Since $m-r\geq1$, one has $h\geq2\epsilon+1>0$.
Direct expansion gives
\begin{equation}\label{eq:pfaffian-critical-exponents}
 -a_k(-\rho)=(p-k)(h+2k)\quad(0\leq k<r),
 \qquad a_k(-\rho)>0\quad(k\geq r).
\end{equation}
In particular, $a_p(-\rho)=0$ and $a_0(-\rho)=-ph$.
The critical slopes are $a_k'(s)=r-k$, for $0\leq k<r$.

For clarity, the continuation needed here is elementary Taylor
subtraction.  If $A(s)$ is holomorphic and $H(t,s)$ is smooth on
$[0,1]$ and holomorphic in $s$, write
$H(t,s)=\sum_{\ell<M}H_\ell(s)t^\ell+t^MR_M(t,s)$.
Then
\[
 \int_0^1t^{A(s)-1}H(t,s)\,dt
 =\sum_{\ell<M}\frac{H_\ell(s)}{A(s)+\ell}
      +\int_0^1t^{A(s)+M-1}R_M(t,s)\,dt.
\]
The last term is holomorphic when $\operatorname{Re}A(s)+M>0$.
Apply this successively to $t_2,\ldots,t_r$ in
\eqref{eq:pfaffian-mellin-coordinates}, taking $M$ larger than all
critical nonnegative Taylor degrees in
\eqref{eq:pfaffian-critical-exponents}.
All derivatives of $L^{-a_0(s)}V(t)U(t)^s$ are uniformly bounded on
the closed cube for $s$ near $-\rho$; the tail exponents remain positive.
More explicitly, apply the Taylor projection $T_j$ through the critical
Taylor degree in each variable and put $R_j=1-T_j$.
The commuting expansion $\prod_{j=2}^r(T_j+R_j)$ expresses the integral
as a finite sum.  Every occurrence of $R_j$ retains an extra vanishing
power of $t_j$, making that coordinate integral holomorphic near $-\rho$;
its mixed derivatives satisfy the same estimate.  The tail exponents
have no $s$-dependence and are positive.
Thus a term involving a remainder has at most $r-1$ pole factors.
There are exactly $r$ possible simple-pole factors in the term with no
remainder, including the Gamma factor.  Consequently the pole order is
at most $r$.

To compute its highest coefficient, every critical residue must be taken.
The residue in $t_r$ uses Taylor degree zero and hence sets $t_r=0$.
Write
\begin{align*}
 L_p&=1+t_2+\cdots+t_2\cdots t_p,&
 V_p&=\prod_{1\leq i<j\leq p}(1-t_{i+1}\cdots t_j)^4,\\
 L_{\rm tail}&=1+t_{r+1}+\cdots+t_{r+1}\cdots t_m,&
 V_{\rm tail}&=\prod_{r\leq i<j\leq m}(1-t_{i+1}\cdots t_j)^4.
\end{align*}
At $t_r=0$ one has
$V=V_pV_{\rm tail}$, $L=L_p$, and $U=L_{\rm tail}$.
Indeed, every Vandermonde factor crossing the cut becomes one, and a term
of $U$ has zero $t_r$-degree exactly when it selects all $p$ head indices
and one tail index.
The remaining head coefficient is
\[
 K_{p,h}=
 \left[\prod_{k=1}^{p-1}t_{k+1}^{(p-k)(h+2k)}\right]V_pL_p^{ph}.
\]
Introducing $z_i/z_1=t_2\cdots t_i$, the identity
$\Delta_p(z)^4=(\prod_i z_i^{4(p-i)})V_p$
converts the selected ratio degrees into equal exponents
$H=h+2(p-1)$.  More explicitly, the degree required in the suffix
$z_{k+1},\ldots,z_p$ is
\[
 4\sum_{i=k+1}^p(p-i)+(p-k)(h+2k)=(p-k)H;
\]
successive differences force every individual degree to be $H$.
It follows that $K_{p,h}$ is the coefficient in
Lemma~\ref{lem:fourth-vandermonde} and is strictly positive.

The tail integral equals
\[
 J_{\rm tail}=
 \int_{(0,1)^{m-r}}
 \prod_{j=r+1}^m t_j^{a_{j-1}(-\rho)-1}
       V_{\rm tail}L_{\rm tail}^{-\rho}\,dt_{r+1}\cdots dt_m.
\]
Its exponent parameters are positive; its integrand is positive in the
open cube and bounded by an integrable product of powers.
Thus $0<J_{\rm tail}<\infty$.
The Gamma residue in the variable $s$ is
$(-1)^{ph}/((ph)!r)$, and the product of the other critical slopes is
$p!$.  Consequently
\[
 \lim_{s\to-\rho}(s+\rho)^rZ_r(s)
 =\frac{c_N(-1)^{ph}J_{\rm tail}}{r!}
       \prod_{j=0}^{p-1}\frac{(2j+2)!}{2!(h+2j)!}\ne0.
\]
This proves the exact pole order.
\end{proof}

\begin{remark}[Compact localization of the highest pole]
\label{rem:compact-mellin}
In Proposition~\ref{prop:pfaffian-mellin}, let
$\chi\in C_c^\infty([0,\infty))$ equal one on a neighborhood of zero,
and set
\[
 Z_{r,\chi}(s)=\int_{\C^{D_0}}Q_r(X)^s\phi_r(X)
                         \chi(\|X\|^2)\,dX.
\]
Then $Z_{r,\chi}$ has the same coefficient of $(s+\rho)^{-r}$
as $Z_r$. In particular, this coefficient is detected by a smooth
test function supported in an arbitrarily small neighborhood of the origin.
\end{remark}
\begin{proof}
Choose $\delta>0$ such that $\chi=1$ on $[0,\delta]$.
Write $t=(t_2,\ldots,t_m)$ and use the ratio coordinates of the proposition.
In the initial convergence region,
\[
 Z_r(s)-Z_{r,\chi}(s)
 =c_N\int_{(0,1)^{m-1}}
       \prod_{j=2}^m t_j^{a_{j-1}(s)-1}H(t,s)\,dt,
\]
where
\[
 H(t,s)=V(t)U(t)^s
       \int_0^\infty t_1^{a_0(s)-1}
       \bigl(1-\chi(t_1L(t))\bigr)e^{-t_1L(t)}\,dt_1.
\]
Since $1\le L(t)\le m$, the integrand of the inner integral is zero
for $0<t_1\le\delta/m$, uniformly on the closed ratio cube.
For every compact set of $s$-values, every mixed derivative in the
ratio variables, and every finite number of $s$-derivatives, the remaining
integrand is bounded by
\[
 C(1+t_1)^A(1+|\log t_1|)^B e^{-t_1},
 \qquad t_1\ge\delta/m,
\]
for constants depending only on the derivative orders and the compact set.
Indeed, derivatives of $L$ are bounded, derivatives of $\chi$ are bounded
and compactly supported, derivatives of the exponential introduce only
powers of $t_1$, and $U\ge1$ on the closed cube.
Thus $H$ is entire in $s$, smooth on the closed ratio cube, and has
uniformly bounded mixed derivatives there for $s$ in compact sets.

At $s=-\rho$, apply the same commuting Taylor projections in
$t_2,\ldots,t_r$ as in the proposition. The tail exponents are positive
and independent of $s$. Each of these $r-1$ critical ratio variables
can contribute at most a simple pole, and there is now no Gamma factor.
Consequently $Z_r-Z_{r,\chi}$ has pole order at most $r-1$ at $-\rho$.
It follows that
\[
 \lim_{s\to-\rho}(s+\rho)^r Z_{r,\chi}(s)
 =\lim_{s\to-\rho}(s+\rho)^r Z_r(s)\ne0,
\]
as claimed. The cutoff may be chosen with arbitrarily small support.
\end{proof}

\subsection{The Bernstein equation and exact integral multiplicities}\label{subsec:integer-saturation}
\begin{lemma}[Poles and total integral multiplicity]
\label{lem:integer-pole-bound}
Let $I=(f_1,\ldots,f_a)\subset\C[z_1,\ldots,z_d]$ be nonzero, and
let $\phi$ be Schwartz.  If the continuation of
\[
 Z_I(\phi;s)=\int_{\C^d}\Bigl(\sum_j|f_j(z)|^2\Bigr)^s\phi(z)\,dz
\]
has a pole of order $R$ at a negative integer $-\rho$, then
\[
 R\leq\sum_{u\in\Z}\operatorname{mult}_{s=-u}b_I(s).
\]
\end{lemma}

\begin{proof}
For a polynomial $g$ choose a global Bernstein equation
$P(s)g^{s+1}=b_g(s)g^s$.
Its $M$-fold iteration, followed by integration by parts in the
holomorphic variables only, gives
\begin{equation}\label{eq:one-sided-bernstein}
 \left(\prod_{j=0}^{M-1}b_g(s+j)\right)
       \int |g|^{2s}\psi
 =\int g^M|g|^{2s}P_M(s)^{\mathsf t}\psi,
 \qquad P_M(s)=P(s)P(s+1)\cdots P(s+M-1).
\end{equation}
The transpose is the bilinear formal transpose, so that
$\partial_{z_i}^{\mathsf t}=-\partial_{z_i}$ and $s$ is fixed.
One first proves the identity for sufficiently large real part of $s$
with cutoffs, where differentiation across $g=0$ is legitimate;
Schwartz decay permits removal of the cutoffs at infinity.
Only holomorphic derivatives are used, and they do not differentiate
$\overline g^{\,s}$ away from $g=0$.
The right side of \eqref{eq:one-sided-bernstein} is holomorphic whenever
$M+2\operatorname{Re}s>0$.
Near $g=0$, its positive power of $|g|$ controls all parameter derivatives,
which add powers of $\log|g|$; at infinity the polynomial coefficients of
$P_M$ are controlled by Schwartz decay.
These observations justify the continuation and the pole-order bound
by the zero order of the single shifted product.

Discard any zero generators from the chosen list, which does not change $I$ or $Z_I$. Now put $g(z,y)=\sum_{j=1}^a y_jf_j(z)$ and
$\psi(z,y)=\phi(z)e^{-\|y\|^2}$.
Musta\c{t}\u{a}'s comparison theorem
\cite[Theorem~1.1]{Mustata} gives
\[
 b_g(s)=(s+1)b_I(s).
\]
A unitary change of variables in $y$, for every $z$ outside the common
zero set of the $f_j$, and elementary one-variable Gaussian integration
give
\begin{equation}\label{eq:gaussian-suspension}
 \int |g|^{2s}\psi
       =\pi^a\Gamma(s+1)Z_I(\phi;s).
\end{equation}
The exceptional set has measure zero.  Fubini proves the formula in its
initial region, and continuation proves it meromorphically.
At $s=-\rho$ the left side has pole order $R+1$.
Choose an integer $M>2\rho$.
Exactly one of the additional factors $s+j+1$, $0\leq j<M$, vanishes
there.  Equations \eqref{eq:one-sided-bernstein} and
\eqref{eq:gaussian-suspension} therefore give
\[
 R+1\leq1+\sum_{j=0}^{M-1}
          \operatorname{mult}_{s=-(\rho-j)}b_I(s)
       \leq1+\sum_{u\in\Z}\operatorname{mult}_{s=-u}b_I(s).
\]
Subtracting one proves the assertion.
\end{proof}

\begin{corollary}\label{cor:integer-total-lower}
For $1\leq r<m$,
\[
 \sum_{u\in\Z}\operatorname{mult}_{s=-u}b_{\pr_r}(s)\geq r.
\]
\end{corollary}
\begin{proof}
For $r\geq2$ combine Proposition~\ref{prop:pfaffian-mellin} with
Lemma~\ref{lem:integer-pole-bound}.
For $r=1$, $Q_1(X)=\|X\|^2$ and direct radial integration gives
\[
 \int_{\C^{D_0}}Q_1(X)^s e^{-\|X\|^2}\,dX
      =\pi^{D_0}\frac{\Gamma(s+D_0)}{\Gamma(D_0)}.
\]
Its simple pole at $s=-D_0$ proves the claim by the same lemma.
\end{proof}

\begin{corollary}[Exact integral multiplicities]\label{cor:integer-saturation}
For $r<m$, Proposition~\ref{prop:upper} gives
\[
 b_{\pr_r}(s)\mid
 B_{N,r}(s)=\prod_{i=0}^{r-1}\prod_{u=0}^{r-i-1}
 \left(s+\frac{\binom{N-2i}{2}+u}{r-i}\right),\qquad r<m.
\]
Then for every integer $q$,
\[
 \operatorname{mult}_{s=-q}b_{\pr_r}(s)
 =\#\left\{0\leq i<r:
       \left\lceil\frac{\binom{N-2i}{2}}{r-i}\right\rceil=q\right\}.
\]
\end{corollary}
\begin{proof}
For each fixed $i$, the $r-i$ consecutive numerator values contain
exactly one multiple of $r-i$.  Thus $B_{N,r}$ has exactly $r$ integral
factors, counted with multiplicity, and the displayed cardinality is its
multiplicity at $-q$.
Divisibility bounds each multiplicity of $b_{\pr_r}$ by that of $B_{N,r}$.
Corollary~\ref{cor:integer-total-lower} forces equality in their total
integral multiplicities.  The nonnegative individual differences must
therefore all vanish.
\end{proof}

\begin{remark}
The construction in Proposition~\ref{prop:pfaffian-mellin} is restricted
to $r<m$; when $r=m$, its parameter $h$ is negative and the proof does
not apply.  The extremal cases must be handled separately.
\end{remark}

\section{Multiplicity of the nonintegral roots}\label{sec:central}
We show that collisions between nonintegral labels do not produce
repeated roots. The central-character facts used below were recalled
in Sections~\ref{subsec:representations}--\ref{subsec:central-support}.

\subsection{The common kernel of the commuting operators}\label{subsec:common-kernel}

Fix a nonintegral $\rho$ for which $-\rho$ is a root of
$b_{\mathfrak p_r}(s)$. Proposition~\ref{prop:upper} gives a primary
decomposition for $s$ on the graph quotient. Project its cyclic generator
to the component where $\vartheta=s+\rho$ is nilpotent and take its Fourier
transform, denoted by $v$. A fixed power of $\vartheta$ annihilates the
whole component. Set
\[
 q=-rs-D_0,\qquad z=s+2N-2,\qquad
 q_0=r\rho-D_0,\qquad z_0=2N-2-\rho\notin\Z.
\]
The relations used in Proposition~\ref{prop:upper} remain valid after this
projection. Here $B^k v$ denotes any word of length $k$ in the commuting
operators $B_{ab}$, and not their sum.

Direct calculation from \eqref{eq:skew-action-commutators} gives
\begin{equation}\label{eq:B-dual-action}
 [A_{uv},B_{ab}]=-\delta_{ua}B_{vb}-\delta_{ub}B_{av}.
\end{equation}
In particular, the formal span of the $B_{ab}$ is the
$\mathfrak g$-representation $T=\bigwedge^2V^*$.
Define
\[
 W_k=U(\mathfrak g)\C[s]\operatorname{span}\{B^k v\},
 \qquad W=\sum_{k=0}^{r-1}W_k.
\]
Commuting the $B_{ab}$ past $U(\mathfrak g)$ with
\eqref{eq:B-dual-action} preserves word length. Thus
$B_{ab}W_k\subset W_{k+1}$, with $W_r=0$, and every length-$r$ word
annihilates $W$. Moreover,
\begin{equation}\label{eq:W-degree}
 W=\bigoplus_{k=0}^{r-1}W_k,
 \qquad E|_{W_k}=q-k.
\end{equation}
Indeed, $E+rs+D_0$ acts on $W_k$ by the distinct scalar $-k$.
The span of the derivative Pfaffians of each fixed size is stable under
commutators with $\mathfrak g$. Consequently the graph relations extend
from $B^kv$ to $W_k$ and give
\begin{equation}\label{eq:W-pfaffian-relation}
 \mathfrak p_{r-k}(\partial)W_k=0,
 \qquad C_{r-k}W_k=0.
\end{equation}

Put $F_0=\bigcap_{a<b}\ker(B_{ab}:W\to W)$. Its decomposition by the
$W_k$ is direct, because $B_{ab}$ raises $k$ by one. If
$w\in F_0\cap W_k$, Lemma~\ref{lem:capelli} and
\eqref{eq:W-pfaffian-relation} imply
\begin{equation}\label{eq:kernel-product}
 \prod_{i=0}^{r-k-1}\bigl(q-k+iz-i(2i-1)\bigr)w=0.
\end{equation}
At $s=-\rho$, two different factors could vanish only if
$z_0=2(i+j)-1$ for distinct integers $i,j$. This contradicts
$z_0\notin\Z$. A nonvanishing factor is invertible on the primary
component, by the finite geometric-series inverse for a scalar plus a
nilpotent operator. A vanishing factor is a nonzero scalar multiple of
$\vartheta$, since its $s$-slope is $-(r-i)$. We obtain
\begin{equation}\label{eq:epsilon-kernel}
 \vartheta F_0=0.
\end{equation}

If $F_0\cap W_k\ne0$, there is a unique vanishing index $i$ in
\eqref{eq:kernel-product}. Put
\[
 a_i=2i-1-z_0,\qquad \lambda^{(i)}=\lambda(i,a_i).
\]
Then $q_0-k=ia_i$. The same Capelli recurrence, now for every
$0\le p\le m$, yields
\begin{align*}
 p!C_pw
 &=\prod_{h=0}^{p-1}\bigl(ia_i+hz_0-h(2h-1)\bigr)w\\
 &=\prod_{h=0}^{p-1}(i-h)(a_i+2h)w
 =p!c_p(i,a_i)w.
\end{align*}
Lemma~\ref{lem:rectangle-character} shows that $F_0\cap W_k$ has the
exact full central character $\chi_{\lambda^{(i)}}$. Thus $F_0$ has a
finite-codimensional central annihilator, with only these finitely many
characters in its support.

\subsection{Separation of the kernel-filtration layers}\label{subsec:separation}

Define $F_{-1}=0$ and, for $j\ge0$, define
\[
 F_j=\{w\in W:B_{\eta_1}\cdots B_{\eta_{j+1}}w=0
                    \text{ for every word of length }j+1\}.
\]
These are $U(\mathfrak g)$-submodules, $F_{r-1}=W$, and applying all
words of length $j$ induces a $\mathfrak g$-equivariant injection
\begin{equation}\label{eq:kernel-layer-injection}
 F_j/F_{j-1}\lhook\joinrel\longrightarrow
 \operatorname{Hom}_{\C}(\operatorname{Sym}^jT,F_0)
 \simeq\operatorname{Sym}^j(\bigwedge^2V)\otimes F_0.
\end{equation}
Indeed, another $B$ kills every image, the kernel consists precisely of
$F_{j-1}$, and equivariance follows from
\eqref{eq:B-dual-action}. Every weight $\nu$ of
$\operatorname{Sym}^j(\bigwedge^2V)$ has nonnegative integral coordinates
and satisfies $\sum_a\nu_a=2j$. Lemma~\ref{lem:central-tensor} implies
that each layer has a finite-codimensional central annihilator and
support contained in
\[
 \Sigma_j=
 \{\chi_{\lambda^{(i)}+\nu}:0\le i<r,
      \ \nu\text{ a weight of }\operatorname{Sym}^j(\bigwedge^2V)\}.
\]

\begin{lemma}\label{lem:skew-support-separation}
If $j\ne\ell$, then $\Sigma_j\cap\Sigma_\ell=\varnothing$.
\end{lemma}

\begin{proof}
The shifted-coordinate multiset in \eqref{eq:HC-multiset} for
$\lambda^{(i)}+\nu$ has exactly $2i$ coordinates in the nonzero class
$-z_0+\Z$ of $\C/\Z$; the remaining coordinates lie in $\Z$. Equality
of two central characters therefore first forces $i=i'$. Equality of
the trace characters then forces
\[
 2ia_i+2j=2ia_i+2\ell,
\]
so $j=\ell$. This argument applies to every nonintegral $\rho$,
including half-integers.
\end{proof}

\begin{proposition}\label{prop:nonintegral-simple}
Every nonintegral root of $b_{\mathfrak p_r}(s)$ is simple.
\end{proposition}

\begin{proof}
By \eqref{eq:epsilon-kernel} and commutativity,
$\vartheta F_j\subset F_{j-1}$ for all $j$: a length-$j$ word takes
$F_j$ into $F_0$, where $\vartheta$ acts as zero. We prove
$\vartheta F_j=0$ by induction. The case $j=0$ has already been proved.
If $\vartheta F_{j-1}=0$, multiplication by $\vartheta$ induces a
$Z$-linear map
\[
 F_j/F_{j-1}\longrightarrow F_{j-1}.
\]
The source has central support contained in $\Sigma_j$. The target
has a finite-codimensional central annihilator and support contained in
$\bigcup_{\ell<j}\Sigma_\ell$: the product of the central annihilators
of its successive layers annihilates it. These two supports are disjoint
by Lemma~\ref{lem:skew-support-separation}. Their annihilator ideals are
therefore comaximal, and the displayed $Z$-linear map is zero. This proves
the induction and gives $\vartheta W=0$.

In particular, $(s+\rho)v=0$. Inverse Fourier transformation gives
$(s+\rho)\delta_\rho=0$ for the projected graph generator. That generator
is cyclic over the degree-zero graph operator algebra, in which $s$ is
central. Hence $s+\rho$ annihilates the whole primary component. Its
exponent in the minimal polynomial is at most one, and it is exactly one
because $-\rho$ was assumed to be a root.
\end{proof}

\section{Occurrence of every root}\label{sec:roots}
We prove occurrence independently of the multiplicity arguments.
The homogeneous-map criterion of \cite[Theorem~3.6]{CZanalytic} is
recalled with a proof on the full graph quotient, including the local
statement needed for rank slices.

\subsection{A homogeneous-map criterion}\label{subsec:criterion}

\begin{lemma}[Homogeneous immersive-map criterion]\label{lem:immersive}
Let $f_1,\ldots,f_a\in\C[x_1,\ldots,x_d]$ be homogeneous of the same positive degree $e$. Suppose the differential of $F=(f_1,\ldots,f_a)$ has rank $d$ at some $x^0$. For $I=(f_1,\ldots,f_a)$,
\[
 b_I\left(-\frac{d+u}{e}\right)=0\qquad(0\le u<e,\ u\in\Z).
\]
The same values are roots of the local polynomial at the origin, including
its analytic-germ version.
\end{lemma}
\begin{proof}
Reorder the generators so that $F_0=(f_1,\ldots,f_d)$ has invertible Jacobian at $x^0$. Let $\varphi$ be its local holomorphic inverse near $v^0=F_0(x^0)$ and put $J(x)=\det(\partial_i f_j)_{j,i=1}^d$. For every polynomial $h$ the divergence identity is
\begin{equation}\label{eq:divergence}
 \frac{(\partial_i h)(\varphi(v))}{J(\varphi(v))}
 =\sum_{j=1}^d\partial_{v_j}
 \left(\frac{(\partial_i f_j)(\varphi(v))h(\varphi(v))}{J(\varphi(v))}\right).
\end{equation}
Indeed, writing $A_{ji}=\partial_i f_j$, the terms differentiating $h$ give the chain rule, and the other terms cancel because
\[
 \sum_{j,k}(A^{-1})_{kj}\partial_kA_{ji}
 =\sum_{j,k}(A^{-1})_{kj}\partial_iA_{jk}=J^{-1}\partial_iJ.
\]

Fix $u\in\{0,\ldots,e-1\}$. For $\alpha\in\N^d$ define
\[
 L_\alpha(h)=\left.\partial_v^\alpha\frac{h(\varphi(v))}{J(\varphi(v))}\right|_{v=v^0},
 \qquad \lambda_\alpha(h)=L_\alpha(h_{e|\alpha|+u}),
\]
where $h_b$ is the homogeneous part of degree $b$. The bilinear Fischer pairing
\[
 \langle P,h\rangle=P(\partial)h\big|_{x=0}
\]
is nondegenerate on every homogeneous component: $\langle x^\alpha,x^\beta\rangle=\alpha!$ if $\alpha=\beta$, and zero otherwise. Thus $\lambda_\alpha$ is represented by a unique homogeneous polynomial $P_\alpha$ of degree $e|\alpha|+u$. Multiplication by $x_i$ is adjoint to $\partial_i$, and multiplication by $f$ is adjoint to $f(\partial)$. Apply $\partial_v^\alpha$ to~\eqref{eq:divergence} and take the appropriate homogeneous component. The result is
\begin{equation}\label{eq:PDEfamily}
 x_iP_\alpha=\sum_{j=1}^d(\partial_if_j)(\partial)P_{\alpha+e_j}.
\end{equation}
Both sides have degree $e|\alpha|+u+1$. Differentiating by $\partial_i$ and summing over $i$, using the Euler identity for the constant-coefficient polynomials $f_j(\partial)$, gives
\begin{equation}\label{eq:eulerfamily}
 \sum_{j=1}^d f_j(\partial)P_{\alpha+e_j}
 =\left(|\alpha|+\frac{d+u}{e}\right)P_\alpha.
\end{equation}
Extend the family to $\alpha\in\N^a$ by setting $P_\alpha=0$ whenever any of the last $a-d$ coordinates is positive. Equations~\eqref{eq:PDEfamily}--\eqref{eq:eulerfamily} then hold with all $a$ generators in the sums, since those added indices always give zero terms.

The polynomial $P_0$ is nonzero. For $u=0$ test $h=1$. For $u>0$, full rank and $e\ge2$ force $x^0\ne0$; a linear form $\ell$ with $\ell(x^0)\ne0$ gives $\lambda_0(\ell^u)\ne0$.
Define on the free graph module
\[
 \Lambda(h\partial_t^\alpha\delta)=\langle P_\alpha,h\rangle.
\]
The graph actions and the two identities imply, for every $w\in M$,
\begin{equation}\label{eq:functional}
 \Lambda(\partial_{x_i}w)=0,\qquad
 \Lambda(sw)=-\frac{d+u}{e}\Lambda(w).
\end{equation}
Moreover $\Lambda(U^1)=0$. A spanning Weyl monomial of $V$-degree at least one can be ordered as
\[
 \partial_x^\gamma x^\beta\partial_t^\theta t^\eta,
 \qquad |\eta|\ge|\theta|+1.
\]
When $\gamma\ne0$, its action is killed by the first part of~\eqref{eq:functional}. Otherwise its action on $\delta$ is $x^\beta f^\eta\partial_t^\theta\delta$, whose polynomial degree is at least $e(|\theta|+1)>e|\theta|+u$, so the Fischer pairing with $P_\theta$ vanishes.
Choose $h_0$ with $\Lambda(h_0\delta)\ne0$. Since $b_I(s)h_0\delta\in U^1$, applying $\Lambda$ proves the claimed root.

For the local assertion, extend the pairing to a convergent power series
$h$ by taking only its homogeneous Taylor part of degree $e|\alpha|+u$.
The identities~\eqref{eq:functional} remain valid, because each is an
identity on finite Taylor jets. A local analytic $V^1$ operator is a finite
sum of terms that can be ordered as
$\partial_x^\gamma\partial_t^\theta t^\eta c(x,t)$, with
$|\eta|\ge|\theta|+1$ and $c$ holomorphic: this follows by factoring the
corresponding power of the finitely generated ideal $(t_1,\ldots,t_a)$
in each coefficient. On $\delta$ its action is
$\partial_x^\gamma\bigl(c(x,f(x))f^\eta\partial_t^\theta\delta\bigr)$.
When $\gamma\ne0$ the derivative identity kills it. Otherwise its
coefficient has Taylor order at least $e(|\theta|+1)$, so the pairing
vanishes. Thus $\Lambda$ also annihilates the analytic $U^1$, while the
same polynomial witness $h_0$ remains nonzero. Applying a local
Bernstein equation proves the local assertion. The full-rank point
$x^0$ was used only to construct the polynomial family $P_\alpha$;
it need not be the origin.
\end{proof}

\subsection{The Pfaffian map and transverse rank slices}\label{subsec:slices}
\begin{lemma}\label{lem:pfimmersion}
For $r<m$, the map $\bigwedge^2\C^N\longrightarrow\bigwedge^{2r}\C^N$ given by $\omega\mapsto\omega^r/r!$ has injective differential at a form of maximal rank.
\end{lemma}
\begin{proof}
Choose a basis $e_1,f_1,\ldots,e_m,f_m$, together with $z$ if $\epsilon=1$, and set $\omega_0=\sum_{j=1}^m e_j\wedge f_j$. Its differential sends $\eta$ to $\omega_0^{r-1}\wedge\eta/(r-1)!$. Write $\eta$ as the sum of terms $a_je_j\wedge f_j$, terms connecting two different pairs, and, in odd size, terms involving $z$ and one paired basis vector.

Any connecting term is detected individually in the exterior coefficient obtained by adjoining $r-1$ whole pairs avoiding its two pairs. Such pairs exist because $r-1\le m-2$. No different term of $\eta$ can contribute to that coefficient: the whole pairs must be supplied by $\omega_0$, and the two unmatched basis vectors by the chosen term. The odd-size terms are similarly detected by adjoining whole pairs avoiding their one paired index. Thus all these coefficients vanish in a kernel vector. The coefficients on unions of $r$ whole pairs then give $\sum_{j\in A}a_j=0$ for every $r$-subset $A$. Comparing subsets differing by one index, possible since $r<m$, makes all $a_j$ equal. Their sums give $ra_j=0$. Hence $\eta=0$.
\end{proof}

Local polynomials below mean algebraic local polynomials: the graph construction is performed over the local ring of the smooth affine variety at the indicated point. We use three basic properties \cite{BMS,LRWW}: invariance under regular coordinate changes and replacement of ideal generators; invariance under adding a smooth factor; and the least-common-multiple description of the global polynomial in terms of its local polynomials. For a homogeneous ideal the local polynomial at the origin equals the global one. Indeed, a Bernstein equation over the local ring at the origin extends to some Zariski open neighborhood $U$ of the origin. For any point $x$, the inverse image of $U$ under $\lambda\mapsto\lambda x$ is a Zariski open subset of $\mathbb A^1$ containing zero, so it contains some $\lambda\ne0$. Dilation by this $\lambda$ preserves local ideal polynomials, since it rescales all generators by nonzero constants. Every local polynomial therefore divides the origin polynomial; the reverse divisibility follows from the least-common-multiple property. 

\begin{lemma}[Transverse rank slice]\label{lem:slice}
At a matrix of rank $2i<2r$, the algebraic germ of $\pr_r^{(N)}$ is, up to regular coordinates and a smooth factor, the germ at zero of $\pr_{r-i}^{(N-2i)}$.
\end{lemma}
\begin{proof}
Choose an invertible $2i\times2i$ principal skew block and write
\[
 X=\begin{pmatrix} A&B\\-B^t&C\end{pmatrix},\qquad S=C+B^tA^{-1}B.
\]
Block-unipotent congruence takes $X$ to $\operatorname{diag}(A,S)$. The substitution $C\mapsto S$ is an invertible regular coordinate change where $\det A\ne0$, and $S=0$ at the chosen point. Congruence preserves Pfaffian ideals by the exterior-power transformation law. In a $2r$-Pfaffian of the block diagonal matrix, the $S$-part has half-size at least $r-i$; larger Pfaffians belong to the ideal of smaller ones by exterior Laplace expansion. Conversely, adjoining all indices of $A$ to any $2(r-i)$-Pfaffian of $S$ yields $\Pf(A)$ times that Pfaffian. Since $\Pf(A)$ is a unit, the two ideals agree. The $A,B$ coordinates supply the smooth factor.
\end{proof}

\begin{proposition}\label{prop:roots}
If $r<m$, every number in $-R_{N,r}$ is a root of $b_{\pr_r}$.
\end{proposition}
\begin{proof}
Lemma~\ref{lem:pfimmersion} and Lemma~\ref{lem:immersive} give the origin window $-(D_0+u)/r$ for $0\le u<r$. For each rank $2i$, $0\le i<r$, Lemma~\ref{lem:slice} reduces to the degree-$(r-i)$ Pfaffian map in $\binom{N-2i}{2}$ variables. It remains nonmaximal, because $r-i<\lfloor(N-2i)/2\rfloor$. The same two lemmas give all values in its window already at the slice origin, by the local part of Lemma~\ref{lem:immersive}. This local polynomial divides $b_{\pr_r}$, proving the proposition.
\end{proof}

\subsection{The maximal cases}\label{subsec:boundary}
\begin{proposition}[The even maximal boundary]\label{prop:even-boundary}
For the full Pfaffian $f_m$ of a generic $2m\times2m$ skew matrix,
\[
 b_{f_m}(s)=\prod_{j=1}^{m}(s+2j-1).
\]
\end{proposition}
\begin{proof}
We supply a lower bound to complement Remark~\ref{rem:maxupper}.
Give every independent entry the standard complex Gaussian law
$\pi^{-1}e^{-|x|^2}\,dx$, and put $G_m(s)=\mathbb E|f_m(X)|^{2s}$
initially for $\Re s>0$. Write
\[
 X=\begin{pmatrix}0&z\\-z^t&W\end{pmatrix},\qquad z\in\C^{2m-1}.
\]
Condition on $z\ne0$ and choose a measurable unitary matrix $U(z)$ with
$zU(z)^t=(\|z\|,0,\ldots,0)$. Such a choice follows by completing
$\overline z/\|z\|$ to the rows of an orthonormal basis, on measurable
Gram--Schmidt patches. Congruence by $\operatorname{diag}(1,U(z))$
does not change the absolute value of the Pfaffian. Conditional on $z$,
the matrix $W'=U(z)WU(z)^t$ has the same standard Gaussian skew law
as $W$, because the exterior-square representation of a unitary matrix
is unitary. This conditional law is independent of $z$.
Expansion by the first row gives absolute value $\|z\|\,|f_{m-1}(B)|$,
where $B$ is the bottom $(2m-2)$-square block of $W'$. Its entries are
independent standard Gaussians; the unused entries integrate to one.
Since $\|z\|^2$ has Gamma distribution of shape $2m-1$, Fubini gives
\[
 G_m(s)=\frac{\Gamma(s+2m-1)}{\Gamma(2m-1)}G_{m-1}(s),\qquad G_0(s)=1.
\]
Thus $G_m(s)=\prod_{j=1}^m\Gamma(s+2j-1)/\Gamma(2j-1)$.
At $s=-(2m-1)$ all $m$ Gamma factors have simple poles and nonzero
residues, so this integral has a pole of exactly order $m$.
Lemma~\ref{lem:integer-pole-bound} bounds that order by the sum of
integer root multiplicities of $b_{f_m}$. The sharper upper bound
in Remark~\ref{rem:maxupper} has precisely the $m$ distinct factors
displayed in the proposition. Hence each must occur, with multiplicity one.
\end{proof}

For odd size $N=2m+1$, the equality
\[
 b_{\pr_m}(s)=\prod_{j=0}^{m-1}(s+2j+3)
\]
is \cite[Theorem~3.9]{LRWW}. Together with the preceding proposition,
this proves both cases of \eqref{eq:boundary}.

\subsection{Completion of the main formula}\label{sec:assemble}
\begin{proof}[Proof of Theorem~\ref{thm:main}]
Assume first $r<m$. Proposition~\ref{prop:upper} excludes every root outside $-R_{N,r}$ and bounds each root multiplicity by its label count. Proposition~\ref{prop:roots} proves that each value in this set occurs. Proposition~\ref{prop:nonintegral-simple} makes each nonintegral root simple. Corollary~\ref{cor:integer-saturation} gives the full label multiplicity at every integral root. This proves~\eqref{eq:mainmult}.

In each $i$-window the numerators are $r-i$ consecutive integers, so their unique integral quotient is $\lceil\binom{N-2i}{2}/(r-i)\rceil$. Thus~\eqref{eq:mainpoly} is the same formula. When $r=m$, Remark~\ref{rem:maxupper} gives the boundary upper bound; Proposition~\ref{prop:even-boundary} and the odd Pfaffian theorem of \cite[Theorem 3.9]{LRWW} give equality. For $r=1$ the nonmaximal formula, or either boundary when $m=1$, gives $s+D_0$.
\end{proof}

\section{Multiplicity bounds and strong motivic monodromy}
\label{sec:multiplicity-monodromy}
\subsection{Local Bernstein--Sato polynomials}\label{subsec:local}
\begin{corollary}[Every local polynomial]\label{cor:local}
Let $X^0$ have rank $2i$. If $i\ge r$, the local ideal is the unit ideal and its polynomial is $1$. If $i<r$, its polynomial is
\[
 b_{\pr_r^{(N)},X^0}(s)=b_{\pr_{r-i}^{(N-2i)}}(s),
\]
and hence is given completely, including multiplicities, by Theorem~\ref{thm:main} with $(N,r)$ replaced by $(N-2i,r-i)$.
\end{corollary}
\begin{proof}
Use Lemma~\ref{lem:slice}, invariance under a smooth factor, and equality of the homogeneous origin polynomial with the global one. Outside the zero locus a generator is a unit, giving the first assertion.
\end{proof}

\subsection{An interval formula for integral multiplicities}\label{subsec:interval}
\begin{corollary}[Integer interval formula]\label{cor:interval}
Suppose $r<m$ and put $a=N-2r\ge2$, $C=\binom a2$. For an integer $\rho$ set $T=\rho-(2a-1)$ and
\[
 H_C(T;r)=\#\{d\in\Z:1\le d\le r,\ 2d^2-Td+C\le0\}.
\]
Then
\begin{equation}\label{eq:intervalmult}
 \mult_{s=-\rho}b_{\pr_r}(s)=H_C(T;r)-H_C(T-1;r).
\end{equation}
If $T<0$ or $T^2<8C$, then $H_C(T;r)=0$. Otherwise
\begin{equation}\label{eq:Hformula}
 H_C(T;r)=\max\left\{0,
 \min\left(r,\left\lfloor\frac{T+\sqrt{T^2-8C}}4\right\rfloor\right)
 -\max\left(1,\left\lceil\frac{T-\sqrt{T^2-8C}}4\right\rceil\right)+1\right\}.
\end{equation}
\end{corollary}
\begin{proof}
Set $d=r-i$ in the first product of~\eqref{eq:mainpoly}. Since
\[
 \frac{\binom{N-2i}{2}}{r-i}=2d+2a-1+\frac Cd,
\]
the contributing indices are exactly $T-1<2d+C/d\le T$. Subtracting the two sublevel-set counts proves~\eqref{eq:intervalmult}. Solving $2d^2-Td+C\le0$ gives~\eqref{eq:Hformula} with the interval $[1,r]$ imposed. 
\end{proof}

For example, the formula gives
\[
 b_{\pr_2^{(6)}}(s)=(s+6)(s+15/2)(s+8),
\]
and
\[
 b_{\pr_3^{(10)}}(s)=(s+14)(s+29/2)(s+15)^2(s+46/3)(s+47/3).
\]
The largest root is
\[
 -\min_{0\le i<r}\frac{\binom{N-2i}{2}}{r-i}.
\]
Its negative agrees with the log canonical threshold obtained from the
resolution data in \cite[Section~4.1]{CSXZ}, together with the
Bernstein--Sato characterization in \cite[Theorem~2]{BMS}.

\subsection{The sharp uniform multiplicity bound}\label{sec:maxmultiplicity}
\begin{proof}[Proof of Theorem~\ref{thm:maxmultiplicity}]
By Theorem~\ref{thm:main}, all nonintegral roots are simple.
Put $a=N-2r\ge2$ and $C=\binom a2$. All counting parameters $a,r,T$ below are integers.
By Corollary~\ref{cor:interval}, the multiplicity of the integral root $-\rho$ is
\[
 \mu_{a,r}(T)=\#\{1\le d\le r:\lceil 2d+C/d\rceil=T\},
 \qquad T=\rho-(2a-1).
\]
We show that $k=\mu_{a,r}(T)$ implies $r\ge k(k+1)/2$.

For $b=2d-a$ one has
\[
 2d+C/d-(2a-1)=\frac{b(b+1)}{2d}\ge0.
\]
Equality occurs at exactly one integral index, $d_0=\lfloor a/2\rfloor$.
Thus $T<2a-1$ gives no indices and $T=2a-1$ gives at most one.
For $T=2a$, the shell is described by $b^2\le a$ with $d_0$ omitted.
Let $n$ be the largest nonnegative integer at most $\sqrt a$ congruent to $a$ modulo two.
The full shell is the integer interval
\[
 \left[\frac{a-n}{2},\frac{a+n}{2}\right]\cap\Z
\]
with $d_0$ deleted; it has $n$ elements.
If its prefix through $r$ has $k=n\ge1$ elements, then
\[
 r\ge\frac{a+n}{2}\ge\frac{n(n+1)}2.
\]
If $1\le k<n$, then $n\ge k+1$ and
\[
 r\ge\frac{a-n}{2}+k-1
 \ge\frac{n(n-1)}2+k-1
 \ge\frac{k(k+1)}2.
\]
This handles the first two shells, including their truncations.

Suppose now that $T\ge2a+1$ and $k\ge2$.
For two distinct indices $d<e\le r$ in a common ceiling shell,
\[
 (e-d)\left|\frac C{de}-2\right|<1.
\]
Consequently
\begin{equation}\label{eq:maxpairbound}
 a(a-1)<6r(r-1),\qquad a<\frac52r.
\end{equation}
For the second assertion, $a\ge5r/2$ would make $a(a-1)>6r(r-1)$.

Set
\[
 D_t=t^2-4a(a-1),\qquad
 \alpha_t=\frac{t-\sqrt{D_t}}4,\qquad
 \beta_t=\frac{t+\sqrt{D_t}}4.
\]
The untruncated shell consists of the integers in
\[
 [\alpha_T,\alpha_{T-1})\ \cup\ (\beta_{T-1},\beta_T].
\]
Its two interval lengths are $L_T/2-1/4$ and $L_T/2+1/4$, where
\[
 L_T=\frac{\sqrt{D_T}-\sqrt{D_{T-1}}}{2}.
\]
Each half-open interval of length $\ell$ contains at most $\lceil\ell\rceil$ integers, so
\begin{equation}\label{eq:maxshellcount}
 k\le
 \left\lceil\frac{L_T}{2}-\frac14\right\rceil+
 \left\lceil\frac{L_T}{2}+\frac14\right\rceil<L_T+2.
\end{equation}
The function $x\mapsto\sqrt{x^2-4a(a-1)}$ is concave on its domain.
Its consecutive increments decrease, whence
\[
 L_T\le L_{2a+1}=\sqrt{2a+\tfrac14}-\sqrt a=:L(a).
\]
Using $a\ge2$ and~\eqref{eq:maxpairbound}, we obtain
\[
 L(a)\le\left(\sqrt{\frac{17}{8}}-1\right)\sqrt a
 <\frac{\sqrt{85}-2\sqrt{10}}4\sqrt r<\frac34\sqrt r.
\]
The last inequality follows by squaring $\sqrt{85}<3+2\sqrt{10}$.
For $r\ge16$, this yields
\[
 k<\frac34\sqrt r+2
 \le\frac{11}{8}\sqrt r-\frac12
 <\sqrt{2r+\frac14}-\frac12.
\]
Since $k$ is integral, the required bound follows.

The smaller values require only two estimates.
The function $L(a)$ is increasing for $a>1/8$.
If $3\le r\le5$, then~\eqref{eq:maxpairbound} gives $a\le11$, and
\[
 L_T\le L(11)=\frac{\sqrt{89}}2-\sqrt{11}<\frac32.
\]
The two interval lengths are less than $1/2$ and $1$, respectively; hence $k\le2$.
If $6\le r\le15$, then $a\le35$, and
\[
 L_T\le L(35)=\frac{\sqrt{281}}2-\sqrt{35}<\frac52.
\]
The interval lengths are less than $1$ and $3/2$, respectively, so $k\le3$, which suffices.
The two displayed inequalities reduce to $9<11$ and $841<875$.
For $r=2$, a repetition would force $C<6$, leaving only $a=2,3$.
The two ceiling values are $(3,5)$ and $(5,6)$, respectively, with no repetition.
The case $r=1$ is immediate. This proves the upper bound in all cases.

Finally, let $\ell=M_r\ge2$ and take $a=\ell^2$, $T=2\ell^2$.
The first-shell description has $n=\ell$ and its last index is
$\ell(\ell+1)/2\le r$, giving exactly $\ell$ labels.
The associated positive root parameter is $\rho=2a-1+T=4\ell^2-1$.
For $M_r=1$, any root attains the bound.
\end{proof}

\subsection{Resolution data}\label{subsec:resolution}
Let $X=\mathbb A^{D_0}$ be the space of alternating matrices and let
$Z_i=\{\operatorname{rank}X\le2i\}$ for $0\le i<r$.
The following is the standard resolution by successive rank strata;
the numerical data are recorded in \cite[Section~4.1, (4.1)]{CSXZ}.

\begin{proposition}\label{prop:resolution}
Successively blow up the strict transforms of $Z_0,\ldots,Z_{r-1}$,
omitting any final blow-up along a divisor. This gives a
principalization $\pi:Y\to X$ whose reduced total transform is a
simple normal crossing divisor with components $E_0,\ldots,E_{r-1}$.
Their numerical data are
\begin{equation}\label{eq:pf-resolution-data}
 \pr_r\mathcal O_Y=\mathcal O_Y\left(-\sum_{i=0}^{r-1}(r-i)E_i\right),
 \qquad
 K_{Y/X}=\sum_{i=0}^{r-1}\left(\binom{N-2i}{2}-1\right)E_i.
\end{equation}
If $N=2r$, the last component is the strict transform of the Pfaffian
hypersurface, with $(N_{r-1},\nu_{r-1})=(1,1)$.
\end{proposition}
\begin{proof}
We recall the local construction to fix the indexing and the
hypersurface boundary. On a chart of the blow-up of the origin,
write $X=tY$ with $y_{12}=1$. The first two rows and columns of $Y$
form an invertible skew block. Skew Schur elimination, as in
Lemma~\ref{lem:slice}, gives regular coordinates in which
\[
 \pr_r(X)=(t^r)\pr_{r-1}(S),
\]
where $S$ is a generic alternating matrix of size $N-2$ and its
coordinates are independent of $t$ and the remaining smooth
coordinates. The strict transforms of the higher rank loci are the
corresponding rank loci for $S$. These charts cover the blow-up:
each chart chooses a nonzero alternating entry as the pivot.
Induction in the matrix size therefore shows that the successive
centers are smooth, their transforms meet the existing boundary
with normal crossings, and the final total ideal is invertible.
For $r=1$ this is the blow-up of the maximal ideal at the origin;
in ambient dimension one that ideal is already a smooth divisor.

At a general point of $Z_i$ of rank exactly $2i$, earlier blow-ups
are isomorphisms. The transverse slice has dimension
$\binom{N-2i}{2}$ and its ideal is generated in degree $r-i$.
Blowing up its origin gives order $r-i$ for the ideal and
discrepancy $\binom{N-2i}{2}-1$. These determine the coefficients
of the divisor dominating $Z_i$. Later blow-ups do not change
these coefficients on its strict transform. In the case $N=2r$,
the last slice has dimension one, so the component is retained as
a strict transform with multiplicity one and discrepancy zero.
This proves \eqref{eq:pf-resolution-data}.
\end{proof}

\subsection{The strong monodromy inequality}\label{subsec:strong}
\begin{lemma}[Coincident resolution slopes]\label{lem:slope-collision}
Put $a=N-2r\ge0$, $C=\binom a2$, and $d=r-i$.
Then
\begin{equation}\label{eq:slope-form}
 \frac{\nu_i}{N_i}=2d+2a-1+\frac Cd,\qquad 1\le d\le r.
\end{equation}
If two distinct indices $d,e$ give the same slope, then
$C=2de$, and their common slope is
$2a-1+2(d+e)$, an odd integer. There are at most two indices
with any fixed slope.
\end{lemma}
\begin{proof}
Expand $\binom{a+2d}{2}/d$ to get \eqref{eq:slope-form}.
Equality of the values at $d$ and $e$ gives
$(d-e)(2-C/(de))=0$. For $d\ne e$ this is $C=2de$,
which gives the displayed common value. Finally the indices
for slope $\rho$ satisfy the quadratic equation
$2d^2-(\rho-2a+1)d+C=0$.
\end{proof}

\begin{proof}[Proof of Theorem~\ref{thm:strong}]
Apply the resolution formula \eqref{eq:motivicresolution} to
Proposition~\ref{prop:resolution}. In each summand a divisor
contributes at most one factor
$(1-\mathbb L^{-\nu_i}T^{N_i})^{-1}$.
Multiplication by their product therefore gives a polynomial,
proving \eqref{eq:motivicdenominator}. The same presentation has
at most $p_{N,r}(\rho)$ denominator factors of slope $\rho$ in
each summand. Definition~\ref{def:motivicorder} gives the first
inequality of \eqref{eq:strongineq}.

Suppose first $r<m$. Every index with slope $\rho$ supplies the
distinct label $(i,0)\in J_{N,r}(\rho)$. If $\rho$ is integral,
Theorem~\ref{thm:main} gives
\[
 p_{N,r}(\rho)\le |J_{N,r}(\rho)|
                  =\mult_{s=-\rho}b_{\pr_r}(s).
\]
If $\rho$ is nonintegral, Lemma~\ref{lem:slope-collision} gives
$p_{N,r}(\rho)\le1$; whenever it is one, the corresponding
label is a root by Theorem~\ref{thm:main}, and that root is simple.
If $p_{N,r}(\rho)=0$, the desired inequality is immediate.

For $r=m$, the slopes are the distinct integers
\[
 2(m-i)-1+2\epsilon,\qquad 0\le i<m.
\]
They are precisely the simple root parameters in
\eqref{eq:boundary}, so the second inequality follows in this
case as well.

In the topological resolution formula
\eqref{eq:topologicaldefinition}, a summand has pole order at most
$p_{N,r}(\rho)$ at $s=-\rho$, and addition cannot increase that
order. Every finite pole is therefore canceled by $b_{\pr_r}(s)$.
The resulting rational function has no finite poles and hence is
a polynomial, proving \eqref{eq:topdivisibility}.
\end{proof}

\begin{corollary}\label{cor:double-poles}
Every motivic or topological pole of a generic Pfaffian ideal has
order at most two. A pole of order two can occur only at a negative
odd integer. All poles in the maximal cases have order at most one.
\end{corollary}
\begin{proof}
Combine Lemma~\ref{lem:slope-collision} with the first inequality
of \eqref{eq:strongineq} and its topological counterpart.
The maximal slopes are distinct.
\end{proof}

\begin{remark}
Pole orders and root multiplicities need not agree. For example,
when $(N,r)=(21,6)$, the slope $35$ has two resolution indices,
whereas the root $-35$ has multiplicity three. The third integral
label has $u>0$, so it contributes to the Bernstein--Sato polynomial
without adding a resolution divisor of that slope.
\end{remark}

\appendix
\section{An independent geometric proof of nonintegral simplicity}
\label{sec:normal-cone-independent}

This section gives a second proof of nonintegral simplicity. It uses neither the bordering operators nor central-character separation, and does not assume the root formula. Its geometric statement is stronger: the part of the entire Verdier specialization with eigenvalues different from one has no nilpotent monodromy. We distinguish this assertion from semisimplicity on individual stalk cohomology groups, which would not suffice.

Write $\mathcal X=\Lambda^2\C^N$, $I=\pr_r$, $Z=V(I)$, and $n=\dim\mathcal X$. With deformation parameter $u$, put
\[
 \mathcal T_I=\operatorname{Spec}\C[\mathcal X][u,I/u],
 \qquad C_I\mathcal X=\{u=0\}\subset\mathcal T_I.
\]
The ring is viewed as a subring of $\C[\mathcal X][u,u^{-1}]$; thus $\mathcal T_I$ is reduced and irreducible, and its open part $u\ne0$ is $\mathcal X\times\C^*$. Constructible complexes and their loci are taken on the underlying reduced complex analytic spaces. In the algebraic arguments we retain the scheme structures of the deformation and its fibers; the deformation function $u$, including its multiplicity along the central fiber, is never replaced by its reduction. Let $\operatorname{Sp}_{I,\lambda}$ denote the generalized $\lambda$-eigensummand of the Verdier specialization of $\C_{\mathcal X}[n]$ on $C_I\mathcal X$. Equivalently, it is the corresponding nearby-cycle object of the constant complex on this open part. Shifts of complexes will not affect the monodromy statements below. Write $N_\lambda=\log(\lambda^{-1}T)$ for its nilpotent monodromy logarithm.

\begin{theorem}\label{thm:nc-semisimplicity}
If $1\le r<\lfloor N/2\rfloor$, then
\[
 N_\lambda=0\quad\hbox{on }\operatorname{Sp}_{I,\lambda}
 \qquad(\lambda\ne1).
\]
Consequently every nonintegral root of the ordinary, unshifted BMS polynomial $b_I(s)$ is simple.
\end{theorem}

\subsection{The closure of an exterior-power image}
Fix $V=\C^M$, $m=\lfloor M/2\rfloor$, and $1\le d<m$. Define
\[
 \Phi_d:\Lambda^2V\longrightarrow\Lambda^{2d}V,
 \qquad\Phi_d(\omega)=\omega^d/d!,\qquad
 Y_d=\overline{\Phi_d(\Lambda^2V)}
\]
with its reduced variety structure. For a basis $v_1,\ldots,v_M$, set $e_j=v_{2j-1}\wedge v_{2j}$. Nonzero scalar factors in the following orbit representatives are immaterial.

\begin{lemma}\label{lem:nc-orbits}
Every nonzero point of $Y_d$ is $\GL(V)$-conjugate to
\begin{equation}\label{eq:nc-orbit-form}
 e_1\wedge\cdots\wedge e_p\wedge
 (e_{p+1}+\cdots+e_{p+q})^{d-p},
 \qquad 0\le p<d,\quad d-p\le q,\quad p+q\le m.
\end{equation}
Conversely, all these forms occur in $Y_d$. When $q=d-p$, the form is decomposable.
\end{lemma}
\begin{proof}
We first record an elementary normal form over $R=\C[[t]]$, with fraction field $K=\C((t))$. Any alternating matrix over $K$ is congruent, by an element of $\GL_M(R)$, to
\[
 t^{a_1}e_1+\cdots+t^{a_b}e_b,
 \qquad a_1\le\cdots\le a_b\in\Z,
\]
with a possible zero remainder. Indeed, choose a nonzero entry of minimal valuation, move it to position $(1,2)$, and absorb its unit by a basis change. Every other entry is divisible by this pivot in $R$. Simultaneous row and column operations over $R$ clear its two rows and columns. The remaining Schur-complement entries still have valuation at least that of the pivot, since each correction is a product of two entries divided by the pivot. Induction proves the assertion, including Laurent exponents.

The maximal-rank orbit in $\Lambda^2V$ is dense, and its image is a dense orbit in $Y_d$. For $0\ne\eta_0\in Y_d$, choose an algebraic curve through $\eta_0$ whose general point lies in this orbit. Normalize and complete the curve. After a finite extension of its fraction field, its generic point lifts to some maximal-rank $\omega(t)$ with $\Phi_d(\omega(t))=\eta(t)$. Such an extension exists because a point of the nonempty fiber over the algebraic closure has finitely many coordinates, all algebraic over the original fraction field. A finite Puiseux extension again has a formal uniformizer.

Apply the preceding normal form. The congruence matrix and its inverse are regular and invertible at the special point, so they only change the $\GL(V)$-orbit of the limit. We can therefore analyze
\[
 \Phi_d(\omega(t))=\sum_{|S|=d}t^{\sum_{j\in S}a_j}e_S,
 \qquad a_1\le\cdots\le a_m.
\]
The $e_S$ are linearly independent. Regularity and a nonzero limit say exactly that $a_1+\cdots+a_d=0$. Put $c=a_d$, let $p$ be the number of exponents strictly below $c$, and let $q$ be the number equal to $c$. Every minimizing $d$-subset contains the first $p$ indices and chooses its remaining $d-p$ indices from the next $q$. The limit is therefore~\eqref{eq:nc-orbit-form}.

Conversely, for $p=0$ choose $q$ exponents equal to zero and all others equal to one. For $p>0$, choose exponents $-(d-p)$ on the first $p$ blocks, $p$ on the next $q$, and $p+1$ on the remaining blocks. The least $d$-sum is zero and gives the asserted limit. This also proves that the list is exhaustive without an assumption about diagonal one-parameter subgroup limits.
\end{proof}

For $\eta\in\Lambda^{2d}V$, set
\[
 K_\eta=\{v\in V:v\wedge\eta=0\},\qquad
 U_d=\{\eta\in Y_d:K_\eta=0\},\qquad
 \Omega_d=\{\omega\in\Lambda^2V:\operatorname{rank}\omega>2d\}.
\]

\begin{lemma}\label{lem:nc-cover}
The set $U_d$ is open and equals $\Phi_d(\Omega_d)$. The map
\[
 \Phi_d:\Omega_d\longrightarrow U_d
\]
is a finite \'etale $\mu_d$-torsor, where $\mu_d$ acts by scalar multiplication. In particular $U_d$ is smooth.
\end{lemma}
\begin{proof}
In~\eqref{eq:nc-orbit-form}, let $A$ be the forced $2p$-dimensional subspace and let $B$ be the symplectic $2q$-dimensional subspace. If $d-p<q$, the exterior multiplication map
\[
 B\longrightarrow\Lambda^{2(d-p)+1}B,
 \qquad v\longmapsto v\wedge\omega_B^{d-p}
\]
is injective, as expansion in a symplectic basis shows. Decomposing exterior powers according to $V=A\oplus B\oplus C$ gives $K_\eta=A$. If $d-p=q$, the form is decomposable and $\dim K_\eta=2d$. Thus $K_\eta=0$ holds exactly for the types $p=0,q>d$, proving $U_d=\Phi_d(\Omega_d)$. It is open because injectivity of the matrix $v\mapsto v\wedge\eta$ is an open condition.

At $\omega\in\Omega_d$ of rank $2q$, decompose $V=S\oplus C$, where $S$ is its support. The derivative
\[
 d\Phi_{d,\omega}(\xi)=\omega^{d-1}\wedge\xi/(d-1)!
\]
is injective. On $\Lambda^2S$ this is the injectivity proved in
Lemma~\ref{lem:pfimmersion}, applied to $\dim S=2q$ and $d<q$.
On $S\otimes C$ it follows from injectivity of
$v\mapsto\omega^{d-1}\wedge v$, which is immediate by expansion
in a symplectic basis; on $\Lambda^2C$ it follows from $\omega^{d-1}\ne0$. The images have different $C$-degrees.

Each fiber is exactly a scalar $\mu_d$-orbit. First, the support of $\omega^d$ equals that of $\omega$ whenever $\operatorname{rank}\omega\ge2d$: contractions of a standard symplectic expression recover every basis vector of the support. For two points of the same fiber, work on their common support $S$. If $H$ stabilizes $\omega^d$, derivative injectivity gives $\operatorname{Lie}H=\mathfrak{sp}(S,\omega)$, hence $H^0=\operatorname{Sp}(S,\omega)$. Every element of $H$ normalizes $H^0$, and the invariant subspace of $\Lambda^2S$ under $\operatorname{Sp}(S,\omega)$ is the line $\C\omega$. Therefore $h\omega=c\omega$ with $c^d=1$. The converse is immediate.

We check properness, rather than inferring it from finite fibers. Suppose $\omega(t)\in\Omega_d(K)$ and $\Phi_d(\omega(t))$ extends to an $R$-point of $U_d$. Apply the Laurent normal form, allowing zero blocks. Its special value lies in $U_d$, so the minimum-sum analysis in Lemma~\ref{lem:nc-orbits} must have $p=0,q>d$. Thus its cutoff exponent is zero and all exponents are nonnegative. The matrix $\omega(t)$ extends over $R$, and its limiting rank is $2q>2d$. This is the valuative properness criterion. Alternatively, in the finite compactification supplied by Zariski's Main Theorem, a boundary point would admit a curve with generic point in $\Omega_d$; the just-proved extension and separatedness would contradict that it is a boundary point. Consequently the quasi-finite map is finite.

The free quotient $\Omega_d/\mu_d$ is smooth, and its map to $U_d$ is finite, bijective on geometric points, and unramified by derivative injectivity. We spell out the scheme-theoretic conclusion used here and below. For a finite unramified map of complex varieties with a single point over every closed point and equal residue fields, each fiber algebra is $\C$: it is finite reduced and has one geometric point. Nakayama's lemma then shows locally that its finite algebra is generated by $1$ over the target ring, so the map is a closed immersion. If the target is reduced and the underlying image is all of it, its kernel is zero. It is therefore an isomorphism. Applying this to $\Omega_d/\mu_d\to U_d$ proves the assertion.
\end{proof}

\subsection{The whole deformation on the good open set}
Let $I_d$ be the ideal generated by the coordinates of $\Phi_d$, and write
\[
 \mathcal T_d=\operatorname{Spec}\C[x,u,I_d/u]
 \subset\Lambda^2V_x\times\Lambda^{2d}V_y\times\mathbf A^1_u.
\]
This is the closure of $y=\Phi_d(x)/u$ for $u\ne0$, with its reduced integral structure. Its projection to the $y$-space factors through $Y_d$. Its normal-cone fiber over $x=0$ is precisely $Y_d$: since the generators all have degree $d$, taking the degree-$kd$ component identifies
\[
 \bigoplus_{k\ge0}I_d^k/\mathfrak m I_d^k
 \simeq\C[\Phi_{d,1},\ldots,\Phi_{d,a}].
\]

\begin{lemma}\label{lem:nc-whole-model}
Over $y\in U_d$ there is an isomorphism
\begin{equation}\label{eq:nc-whole-model}
 \mathcal T_d|_{U_d}\simeq
 (\Omega_d\times\mathbf A^1_v)/\mu_d,
 \qquad (z,v)\longmapsto(x=vz,\ y=\Phi_d(z),\ u=v^d),
\end{equation}
where $\zeta\in\mu_d$ acts by $(z,v)\mapsto(\zeta z,\zeta^{-1}v)$. In a local analytic trivialization of the torsor, the deformation function is $u=v^d$ on a smooth product.
\end{lemma}
\begin{proof}
The action is free and the source quotient is smooth. The displayed map lands in the deformation, since it does so for $v\ne0$, a dense open. The quotient is finite over $U_d\times\mathbf A^1_u$, using Lemma~\ref{lem:nc-cover} and $u=v^d$. Its map to $\mathcal T_d|_{U_d}$ is also finite: on affine opens, generators as a finite module over the base ring remain generators over the intervening target ring.

For $u\ne0$, choosing $v^d=u$ and setting $z=x/v$ gives the inverse, independently of the choice after taking the quotient. The finite image is closed and contains this dense open of the reduced target, so it is the entire target. Over $u=0$, every geometric point of the source has $v=0$ and $x=0$; the torsor fiber of $z$ becomes one point after quotienting. Hence every geometric fiber of the finite map consists of a single point. This also excludes hidden central points with $x\ne0$ on this target open.

At $v=0$ its derivative into the ambient space is
\[
 (\dot z,\dot v)\longmapsto
 (z\dot v,\ d\Phi_{d,z}(\dot z),\ 0)
\]
if $d>1$; for $d=1$ the last entry is $\dot v$. Since $z\ne0$ and $d\Phi_{d,z}$ is injective, this derivative is injective. The quotient map is \'etale, so the induced map is unramified as well. Away from $v=0$ it is already an isomorphism. The finite-map argument at the end of Lemma~\ref{lem:nc-cover}, with the reduced target, now proves~\eqref{eq:nc-whole-model}.
\end{proof}

\subsection{Bad points and their monodromy}
The group $\GL(V)\times\C^*$ acts on $\mathcal T_d$ by
\begin{equation}\label{eq:nc-action}
 (g,t)\cdot(x,y,u)=(gx,t^{-1}gy,tu).
\end{equation}
If $0\ne\eta\in Y_d\setminus U_d$ is not decomposable, Lemma~\ref{lem:nc-orbits} gives it a nonzero forced volume factor. Scaling just one basis vector in this factor defines an algebraic $g(t)$ with $g(t)\eta=t\eta$. If $\eta$ is decomposable, use any basis vector in its support. Combined with the second factor in~\eqref{eq:nc-action}, this fixes $(0,\eta,0)$ and multiplies $u$ by $t$. For $\eta=0$, the second factor alone has the same properties.

We use the following elementary local observation. If an algebraic $\C^*$-action fixes a point $p$ of the central fiber of a complex family $h:\mathcal Y\to\C$ and satisfies $h(tq)=t h(q)$, then the nearby-cycle stalk at $p$, with constant coefficients on the punctured family, has identity monodromy. Indeed, embed an invariant affine neighborhood equivariantly in a finite-dimensional representation and choose an invariant Hermitian metric for the compact circle. A small ball centered at $p$ is circle-invariant. The circle action trivializes the family over a small circle in the base and returns to the identity after one turn. This proves the statement on every cohomology group of the local Milnor fiber, also when $\mathcal Y$ is singular. In particular, the entire stalk of each generalized eigensummand with eigenvalue different from one is zero.

We now apply these facts on every rank chart of $\mathcal X$. At a point of rank $2i<2r$, skew Schur elimination gives a local product
\begin{equation}\label{eq:nc-slice}
 (\mathcal X,I)\simeq(\text{smooth factor})\times
 (\Lambda^2\C^{N-2i},I_{r-i}).
\end{equation}
This is an equality of ideals on the chart, so the deformation to the normal cone has the same product decomposition. The nonmaximal hypothesis ensures $r-i<\lfloor(N-2i)/2\rfloor$.

Call a central point good if its normal direction on this slice belongs to $U_{r-i}$. Lemma~\ref{lem:nc-whole-model} gives a whole open neighborhood of each good point, with deformation function $u=v^{r-i}$. Every point of its central fiber has transverse coordinate $x=0$, so throughout the resulting neighborhood in the central fiber the base rank remains $2i$. Thus the union $G$ of all good points is an open subset of the entire normal cone, with the reduced analytic structure used for constructible complexes; no passage between distinct base ranks has been omitted. On $G$, the nearby-cycle monodromy is of finite order locally as an endomorphism of the entire local-system object, and its nilpotent part is zero.

At every point of the complement $C_I\mathcal X\setminus G$, apply the preceding weight-one action on its own slice, extending it by the identity on the smooth factor in~\eqref{eq:nc-slice}. This is an actual local action on the full deformation chart, fixes the point, and has weight one on $u$. Hence the nearby-cycle stalks for eigenvalues different from one all vanish there.

Fix $\lambda\ne1$, write $K=\operatorname{Sp}_{I,\lambda}$, and let $j:G\hookrightarrow C_I\mathcal X$ and $i:C_I\mathcal X\setminus G\hookrightarrow C_I\mathcal X$. We have $i^*K=0$. The localization triangle therefore gives
\begin{equation}\label{eq:nc-clean}
 K\simeq j_!j^*K.
\end{equation}
On the smooth open $G$, the object $j^*K$ is a local system in one cohomological degree, up to the common shift. Its locally vanishing nilpotent endomorphism is therefore zero globally on $G$. Functoriality of~\eqref{eq:nc-clean} now implies $N_\lambda=0$. This proves the geometric assertion of Theorem~\ref{thm:nc-semisimplicity}. Notice why it proves more than pointwise semisimplicity: the whole $\lambda\ne1$ object disappears on the bad boundary, so extension by zero rules out endomorphisms invisible on stalk cohomology.

\subsection{The comparison with the ordinary BMS lattice}
We give the comparison with its index shift, rather than using only equality of eigenvalue sets. Let $a$ be the number of generators of $I$, let $M=(i_f)_+\mathcal O_{\mathcal X}$ on $W=\mathcal X\times\mathbf A^a_t$, and put
\[
 \theta_t=\sum_{j=1}^a t_j\partial_{t_j},\qquad s=-\theta_t-a.
\]
Deform the smooth ambient space $W$ along $H=\mathcal X\times\{0\}$. The closure of the graph family in this deformation is exactly the intrinsic Rees deformation $\mathcal T_I$. Its inclusion is closed and hence proper. Compatibility of nearby cycles with this proper direct image identifies the intrinsic specialization just computed with the canonical graph specialization, including its monodromy endomorphism.

The normalization is as follows. Let $L$ be the pullback of $M$ on the punctured ambient deformation, extended by the open direct image, and let $\rho$ be its projection to $W$. By \cite[\S1.3, (1.3.1)]{BMS},
\begin{equation}\label{eq:nc-v-shift}
 \rho_\bullet\operatorname{gr}_V^\gamma L
 =\bigoplus_{k\in\Z}
 \operatorname{gr}_V^{\gamma+a-1-k}M\otimes u^k.
\end{equation}
On the summand indexed by $k$, put $\beta=\gamma+a-1-k$. Equation~(1.3.2) of \cite{BMS} identifies the deformation Euler operator $\Theta_u$ with $\theta_t+k$, and consequently
\begin{equation}\label{eq:nc-nilpotent-shift}
 \Theta_u-\gamma+1
 =\theta_t+a-\beta=-(s+\beta).
\end{equation}
By the Riemann--Hilbert nearby-cycle comparison \cite[Theorem~3.9]{SaitoIntro}, $T$ corresponds to $\exp(-2\pi i\partial_u u)$. Taking $0<\gamma\le1$, its eigenvalue on this summand is $\lambda=\exp(-2\pi i\gamma)=\exp(-2\pi i\beta)$, and
\[
 N_\lambda=-2\pi i(\Theta_u-\gamma+1)=2\pi i(s+\beta).
\]
Every rational $\beta$ occurs in~\eqref{eq:nc-v-shift} for an appropriate $k$ and such a representative $\gamma$. If $\beta\notin\Z$, its monodromy eigenvalue is different from one. The geometric conclusion therefore yields
\begin{equation}\label{eq:nc-v-annihilation}
 (s+\beta)\operatorname{gr}_V^\beta M=0
 \qquad(\beta\notin\Z).
\end{equation}

Finally let $Q=U^0/U^1$ be the ordinary graph quotient of~\eqref{eq:Qdef}. The canonical filtration induces
\[
 F^\beta Q=((U^0\cap V^\beta M)+U^1)/U^1.
\]
It is finite by \cite[\S2.1, (2.1.2)]{BMS}, and each graded piece is a subquotient of $\operatorname{gr}_V^\beta M$. Fix a nonintegral root $-\rho$ and take the $(s+\rho)$-primary summand $Q_\rho$ of $Q$. The primary projection is a polynomial in $s$ and preserves this filtration. For $\beta\ne\rho$, the operator $s+\beta$ is invertible on $Q_\rho$; its inverse is a finite polynomial in $s+\rho$ and also preserves the filtration. But $s+\beta$ is nilpotent on the $\beta$-graded piece, so that piece vanishes. Thus the finite filtration on $Q_\rho$ has only its exact $\rho$-grade potentially nonzero. Equation~\eqref{eq:nc-v-annihilation} kills $s+\rho$ on that grade, and hence on $Q_\rho$ itself. The root $-\rho$ is simple, completing the theorem.

\begin{remark}
The geometric assertion concerns the entire specialization object.
Its implication for the ordinary Bernstein--Sato polynomial uses the
finite filtration on $Q$. The converse need not hold: simplicity of
the roots in one residue class alone does not imply vanishing of the
nearby-cycle logarithm.
\end{remark}

\phantomsection


\addcontentsline{toc}{section}{References}
\begin{thebibliography}{99}
\small
\bibitem{ADF}
S.~Abeasis and A.~Del Fra,
\emph{Young diagrams and ideals of Pfaffians},
Adv. Math. \textbf{35} (1980), no.~2, 158--178.


\bibitem{BG}
J.~N.~Bernstein and S.~I.~Gelfand,
\emph{Tensor products of finite and infinite dimensional representations of semisimple Lie algebras},
Compositio Math. \textbf{41} (1980), no.~2, 245--285.


\bibitem{BMS}
N.~Budur, M.~Musta\c t\u a, and M.~Saito,
\emph{Bernstein--Sato polynomials of arbitrary varieties},
Compositio Math. \textbf{142} (2006), no.~3, 779--797.


\bibitem{CSS}
S.~Caracciolo, A.~D.~Sokal, and A.~Sportiello,
\emph{Algebraic/combinatorial proofs of Cayley-type identities for derivatives of determinants and pfaffians},
Adv. Appl. Math. \textbf{50} (2013), no.~4, 474--594.
\href{https://arxiv.org/abs/1105.6270}{arXiv:1105.6270}.

\bibitem{CSXZ}
Y.~Chen, Q.~Shi, Y.~Xu, and H.~Zuo,
\emph{On the monodromy conjecture, holomorphy conjecture, and embedded Nash problem for Pfaffian ideals},
preprint (2025).
\href{https://arxiv.org/abs/2511.01262}{arXiv:2511.01262}.

\bibitem{CZanalytic}
Y.~Chen and H.~Zuo,
\emph{Analytic criteria for detecting Bernstein--Sato roots},
preprint (2026).
\href{https://arxiv.org/abs/2610.14083}{arXiv:2610.14083}.

\bibitem{DL}
J.~Denef and F.~Loeser,
\emph{Germs of arcs on singular algebraic varieties and motivic integration},
Invent. Math. \textbf{135} (1999), no.~1, 201--232.
\href{https://arxiv.org/abs/math/9803039}{arXiv:math/9803039}.

\bibitem{Morris}
I.~M.~Gessel, L.~Lv, G.~Xin, and Y.~Zhou,
\emph{A unified elementary approach to the Dyson, Morris, Aomoto, and Forrester constant term identities},
J. Combin. Theory Ser. A \textbf{115} (2008), no.~8, 1417--1435.


\bibitem{Hashimoto}
T.~Hashimoto,
\emph{On the principal symbols of $K_{\C}$-invariant differential operators on Hermitian symmetric spaces},
J. Math. Soc. Japan \textbf{63} (2011), no.~3, 837--869.


\bibitem{HU}
R.~Howe and T.~Umeda,
\emph{The Capelli identity, the double commutant theorem, and multiplicity-free actions},
Math. Ann. \textbf{290} (1991), 565--619.


\bibitem{KW}
K.~Kinoshita and M.~Wakayama,
\emph{Explicit Capelli identities for skew symmetric matrices},
Proc. Edinburgh Math. Soc. (2) \textbf{45} (2002), no.~2, 449--465.


\bibitem{LRWW}
A.~C.~L\H orincz, C.~Raicu, U.~Walther, and J.~Weyman,
\emph{Bernstein--Sato polynomials for maximal minors and sub-maximal Pfaffians},
Adv. Math. \textbf{307} (2017), 224--252.


\bibitem{Mustata}
M.~Musta\c t\u a,
\emph{Bernstein--Sato polynomials for general ideals vs. principal ideals},
Proc. Amer. Math. Soc. \textbf{150} (2022), no.~9, 3655--3662.


\bibitem{Perlman}
M.~Perlman,
\emph{Regularity and cohomology of Pfaffian thickenings},
J. Commut. Algebra \textbf{13} (2021), no.~4, 523--548.


\bibitem{SaitoIntro}
M.~Saito,
\emph{Introduction to a theory of $b$-functions},
RIMS K\^oky\^uroku \textbf{1550} (2007), 89--103.

\end{thebibliography}
\end{document}